\documentclass[11pt]{article}
\usepackage[a4paper,top=28mm,bottom=28mm,inner=30mm,outer=30mm]{geometry}
 
\usepackage[svgnames,x11names]{xcolor}
\usepackage{amsfonts,amsmath,amsthm,boldline,colortbl,comment,enumerate,float,graphicx,import,makecell,mathtools,multirow,pdflscape,subcaption,tikz}
\mathtoolsset{centercolon}
\usetikzlibrary{through,intersections,calc}
\usepackage[hidelinks]{hyperref}
\usepackage{cite}
\usepackage{booktabs,microtype}
\usepackage[T1]{fontenc}
\usepackage{lmodern}
\newif\ifshowchanges
\showchangestrue

\newtheorem{theorem}{Theorem}
\newtheorem{lemma}[theorem]{Lemma}
\newtheorem{corollary}[theorem]{Corollary}
\newtheorem{proposition}[theorem]{Proposition}

\theoremstyle{definition}
\newtheorem{remark}[theorem]{Remark}
\newtheorem{definition}[theorem]{Definition}

\title{Permutation codes via spectral graph theory}

\author{%
Aida Abiad\thanks{%
  \texttt{a.abiad.monge@tue.nl},
  Department of Mathematics and Computer Science,
  Eindhoven University of Technology, The Netherlands;
  Department of Mathematics and Data Science of
  Vrije Universiteit Brussel, Belgium%
}
\and
Delio Jaramillo-Velez\thanks{\texttt{djaramil@ull.edu.es}
  Department of Mathematics, Statistics and Operations Research,
  Universidad de La Laguna,
  San Crist\'obal de La Laguna, Spain%
}
\and
Kirsten D. Morris\thanks{\texttt{kdmorris@vt.edu}
  Department of Mathematics,
  Virginia Tech, Blacksburg, VA, USA%
}
\and
Joachim Rosenthal\thanks{\texttt{rosenthal@math.uzh.ch}
  Institute of Mathematics,
  University of Zurich, Zurich, Switzerland%
}
}

\date{}

\providecommand{\tr}{\operatorname{tr}}

\newcommand{\inv}{\operatorname{inv}}
\newcommand{\Cay}{\operatorname{Cay}}
\newcommand{\Spec}{\operatorname{Spec}}
\newcommand{\sgn}{\operatorname{sgn}}

\newcommand{\one}{\mathbf{1}}

\begin{document}
%\linenumbers
\maketitle
\begin{abstract}
We use the Eigenvalue Method to obtain upper spectral bounds on the size of permutation codes in the Hamming, Ulam, Cayley and Kendall metrics. For the Kendall metric, we improve the double-ball upper bound of Buzaglo and Etzion [\emph{IEEE Transactions on Information Theory}, 2015] at minimum distance four and lengths seven and eight. For the Cayley metric, we improve the sphere packing bound, and determine when these spectral improvements occur. We derive theoretical spectral criteria for improvements from the ratio-type bounds of low degree, and show that the quadratic ratio-type bound for the Kendall metric always equals sphere packing. Furthermore, we prove sharp spectral bounds and recover known optimal code sizes for the Hamming and Ulam metrics.\\

\end{abstract}

%Submitted at ACM on 2026-09-09
%https://ef.msp.org/submit_new.php?j=aims_amc

%%%%%%%%%%%%%%%%%%%%%%%%%%%%%%%%%%%%%%%%%%%%%%%%%%%%%%%%%%%%%%%%%%%%%%%%%%%%%%%%%
\section{Introduction}
%%%%%%%%%%%%%%%%%%%%%%%%%%%%%%%%%%%%%%%%%%%%%%%%%%%%%%%%%%%%%%%%%%%%%%%%%%%%%%%%%

{
\emph{Permutation codes} are subsets of the symmetric group $S_n$ whose distinct elements are separated by a prescribed minimum distance. They model coding with ordered symbols and arise in rank-based representations, communication channels and data storage~\cite{chadwick1969rank,diaconis1977spearman,kendall1990rank,karlof1989permutation,blake1979coding,colbourn2004permutation}. Rank modulation for flash memories particularly motivates metrics based on adjacent transpositions and translocations~\cite{wang2009capacity,barg2010codes,FarnoudSkachekMilenkovic2013}. The choice of metric determines the permitted errors and the combinatorial structure of the code problem.

We consider the Hamming, Ulam, Cayley and Kendall $\tau$ metrics, denoted by $d_H,d_U,d_C$ and $d_K$. For $M\in\{H,U,C,K\}$, let $A_M(n,d)$ be the largest cardinality of a subset of $S_n$ with minimum $d_M$-distance at least $d$. Our aim is to bound these quantities from above using spectral graph theory, and show that this new approach can bring improvements over previously known bounds, and also provide several cases for which the spectral bounds are attained.

Our main improvements over the best known upper bounds concern minimum Kendall distance four. In particular, we prove
\[
 A_K(7,4)\le419,\qquad A_K(8,4)\le2879
\]
in Propositions~\ref{prop:K74} and~\ref{prop:K84}, respectively. These improve the double-ball upper bounds $420$ and $2880$ of Buzaglo and Etzion~\cite[Corollary~6]{BuzagloEtzion2015} (note that these existing bounds are also mentioned as the best known upper bounds in Table~II of Parvaresh et al.~\cite{ParvareshEtAl2025}).

Our starting point is to formulate the minimum distance condition in coding theory in graph-theoretic terms. An \emph{independent set} in a graph $G$ contains no adjacent pair of vertices, and the \emph{independence number} $\alpha(G)$ is the largest size of such a set. A \emph{$k$-independent set} in a graph $G$ has distinct vertices at pairwise shortest path distance greater than $k$; its largest possible size is the \emph{$k$-independence number} $\alpha_k(G)$. The \emph{graph power} $G^k$ joins distinct vertices at distance at most $k$, so $\alpha_k(G)=\alpha(G^k)$. We will use the following two graph formulations:

\begin{enumerate}[(i)]
\item For each of the four metrics, the \emph{forbidden-distance graph} (or equivalently, the \emph{$(d-1)$-graph power}), denoted $\Gamma_M(n,d)$, has vertex set $S_n$ and joins distinct permutations at $d_M$-distance less than $d$. A code is exactly an independent set in this graph, giving
\begin{equation}\label{eq:intro-independent}
 A_M(n,d)=\alpha(\Gamma_M(n,d)).
\end{equation}
\item For $M\in\{U,C,K\}$, the \emph{base graph}, denoted $G_M(n)$, joins permutations at $d_M$-distance one. Its shortest path distance is $d_M$, and hence, for $d\ge2$,
\begin{equation}\label{eq:intro-higher-independent}
 \begin{aligned}
 \Gamma_M(n,d)&=G_M(n)^{d-1},\\
 A_M(n,d)&=\alpha\bigl(G_M(n)^{d-1}\bigr)
 =\alpha_{d-1}(G_M(n)).
 \end{aligned}
\end{equation}
\end{enumerate}
The identities~\eqref{eq:intro-independent} and~\eqref{eq:intro-higher-independent} give the equivalent graph formulations underlying our strategy. They allow us to estimate the maximum code size using graph-theoretic techniques, in particular eigenvalue bounds, applied to the independence number of the forbidden-distance graph or, for $M\in\{U,C,K\}$, to the $(d-1)$-independence number of the base graph. We use both formulations throughout the paper.

We should note that for Hamming distance, the corresponding construction uses the ambient \emph{Hamming graph} $\mathcal H_n=H(n,n)$ on all words in $[n]^n$, where $[n]=\{1,\ldots,n\}$:
\[
 \Gamma_H(n,d)=\bigl(\mathcal H_n^{\,d-1}\bigr)[S_n].
\]
Here the brackets denote the induced subgraph on permutation vertices; the power is formed before this restriction. Distinct permutations cannot differ in a single position, so the distance-one graph on $S_n$ itself is edgeless (Remark~\ref{rem:hamming}). Forbidden-distance graphs are also called \emph{Gilbert graphs}; see~\cite[Section III]{jiang2004asymptotic} for binary words and~\cite[Definition 2]{YangChenYuan2008} for permutations.

With these formulations in place, we apply the \emph{Eigenvalue Method}, which has been recently introduced by Abiad, Peters and Ravagnani~\cite{APR2025}. The method combines polynomial spectral bounds for the graph $k$-independence parameter, and their corresponding optimization methods. Its main tools are the inertia-type and ratio-type bounds, together with \emph{mixed-integer linear programming} (MILP) and \emph{linear programming} (LP), respectively. For a polynomial $p$ of degree at most $k$, the matrix $p(A)$ (where $A$ is the \emph{adjacency matrix} of the graph $G$) has zero off-diagonal entries between vertices at distance greater than $k$, which allows us to bound $\alpha_k(G)$ using its eigenvalues. For the classical independence number, we will use degree one directly on a forbidden-distance graph. For higher independence we use degree at most $d-1$ on a base graph. Other coding-theoretic applications of the spectral method have recently appeared in~\cite{Abiad2023EigenvalueCodes,Abiad2024EigenvalueCodesb,Abiad2024EigenvalueCodes,byrne2026multilinear}.

In this paper we first show the graph properties needed to apply the method~\cite[Section 3]{APR2025}. The higher-independence formulation requires equality between metric distance and base-graph distance. We prove that the Ulam, Cayley and Kendall base graphs have the required properties, and that all four forbidden-distance graph families are vertex-transitive and hence walk-regular. The Hamming forbidden-distance graph is connected for $d\ge3$; its distance-two case is edgeless and gives the trivial exact value $A_H(n,2)=n!$. Thus the nontrivial applications below satisfy the required conditions. We should note that distance-regularity and an association scheme are not needed for the Eigenvalue Method to be used~\cite{APR2025}.

For $M\in\{U,C,K\}$, the equivalent formulations~\eqref{eq:intro-independent} and~\eqref{eq:intro-higher-independent} encode the same code problem, but applying spectral bounds to the different graphs can give different estimates, since in general, the spectral of a graph and its graph power are not related} The adjacency matrix of a graph power is generally not a polynomial in the base adjacency matrix, see e.g.~\cite{AbiadCFNZ2022}. Our Ulam computations include examples favouring each approach (Table~\ref{tab:ulam}). Methods that exploit additional algebraic structure provide complementary comparisons: instances of it are coherent configurations used for Kendall distance~\cite{LimHagiwara2012}, and the conjugacy-class linear programming used for Cayley distance~\cite{Tarnanen1999,Bogaerts2010}.

We next compare the resulting spectral upper bounds with the known bounds for each metric. For the Hamming distance case, we compare with linear programming, semidefinite programming and combinatorial bounds~\cite{Tarnanen1999,Bogaerts2010,BogaertsDukes2014,BierbrauerMetsch2013,Montemanni2014,Barta2022}. For the Ulam distance, we use Singleton bounds and published exact computations~\cite{FarnoudSkachekMilenkovic2013,Riet}. For the Kendall distance, the comparisons include sphere packing, double balls, coherent configurations and representation-theoretic bounds~\cite{BuzagloEtzion2015,LimHagiwara2012,Abdollahi2023PermutationCodes,ParvareshEtAl2025,Vijayakumaran2016}. We also compute the full conjugacy-class linear program for the Cayley metric.

For the Cayley distance case, we derive explicit formulas for the quadratic ratio-type bound and the cubic ratio-type bound (Proposition~\ref{prop:cayleyratio}) and determine all lengths at which the quadratic ratio-type bound improves distance-three sphere packing: $n=6,8,10,12$ (Corollary~\ref{cor:cayley-small-positive}). The improvements are, respectively,
\[
 45\longrightarrow40,\quad1390\longrightarrow1344,\quad
 78886\longrightarrow75600,\quad7149277\longrightarrow7044141.
\]
For each metric, we also make theoretical comparisons with the best known upper bounds at small minimum distances, where closed formulas are available for the corresponding ratio-type bounds. These comparisons give criteria for improvement and identify the limitations of the Eigenvalue Method (Section~\ref{sec:new-bounds}).

We also show that the inertia-type and ratio-type bounds are sharp in several cases. For example, the quadratic inertia-type bounds give $A_C(4,3)=2$ and $A_K(3,3)=2$ (Propositions~\ref{prop:C43} and~\ref{prop:K33}). In the Cayley example the optimal quadratic ratio-type bound is three, so the inertia-type bound is stronger (we note that not so many metric examples where found previously where the inertia-type bound would outperform the ratio-type bound, see \cite{APR2025}). An affine construction attains the Cayley quadratic ratio-type bound $A_C(5,3)=10$ (Proposition~\ref{prop:C5sharp}). For the Hamming distance, the classical ratio bound recovers $A_H(n,3)=n!/2$, $A_H(5,4)=20$ and $A_H(6,4)=120$ (Proposition~\ref{prop:hsharp}). For the Ulam distance, the ratio bound on the corresponding graph power gives $A_U(n,n-1)=2$ for every $n\ge3$ (Proposition~\ref{prop:Usharp}). The base-graph Ulam bounds are sharp after flooring at $A_U(4,2)=6$ and $A_U(4,3)=2$, as is the Kendall cubic bound at $A_K(4,4)=3$.

This paper is organized as follows. Section~\ref{sec:preliminaries} introduces the investigated permutation metrics and also the needed background on spectral graph theory. Section~\ref{sec:known-bounds} states the coding bounds used for comparison. Section~\ref{sec:permutation-graphs} investigates the needed graph properties and derives the specific spectral bounds. Section~\ref{sec:new-bounds} then derives and compares the bounds for the four metrics, and Section~\ref{sec:conclusion} discusses the resulting strengths and limitations.
\par

%%%%%%%%%%%%%%%%%%%%%%%%%%%%%%%%%%%%%%%%%%%%%%%%%%%%%%%%%%%%%%%%%%%%%%%%%%%%%%%%%
\section{Preliminaries}\label{sec:preliminaries}
%%%%%%%%%%%%%%%%%%%%%%%%%%%%%%%%%%%%%%%%%%%%%%%%%%%%%%%%%%%%%%%%%%%%%%%%%%%%%%%%%

%%%%%%%%%%%%%%%%%%%%%%%%%%%%%%%%%%%%%%%%%%%%%%%%%%%%%%%%%%%%%%%%%%%%%%%%%%%%%%%%%
\subsection{Permutation codes under the most common metrics}\label{sec:metrics}
%%%%%%%%%%%%%%%%%%%%%%%%%%%%%%%%%%%%%%%%%%%%%%%%%%%%%%%%%%%%%%%%%%%%%%%%%%%%%%%%%
Let $n\ge2$, write $[n]=\{1,\ldots,n\}$, and let $S_n$ be the {\emph{symmetric group}} on $[n]$, with identity $e$. We write permutations in one-line notation and use the composition convention $(\sigma\tau)(i)=\sigma(\tau(i))$. We first define the four metrics and then introduce the corresponding extremal code sizes. The one-line word for $\sigma$ is $(\sigma(1)\cdots\sigma(n))$. With our convention, right multiplication by $(i\,j)$ swaps the entries in positions $i$ and $j$, whereas left multiplication relabels symbols.

The symmetric group $S_n$ can be equipped with the following metrics:

%%%%%%%%%%%%%%%%%%%%%%%%%%%%%%%%%%%%%%%%%%%%%%%%%%%%%%%%%%%%%%%%%%%%%%%%%%%%%%%%%
\subsubsection{Hamming metric}
%%%%%%%%%%%%%%%%%%%%%%%%%%%%%%%%%%%%%%%%%%%%%%%%%%%%%%%%%%%%%%%%%%%%%%%%%%%%%%%%%
The \emph{Hamming distance} is
\[
 d_H(\sigma,\tau)=|\{i\in[n]:\sigma(i)\ne\tau(i)\}|.
\]
For example, $d_H(12345,13254)=4$. It measures position-wise disagreement. Its smallest nonzero value on $S_n$ is two.

%%%%%%%%%%%%%%%%%%%%%%%%%%%%%%%%%%%%%%%%%%%%%%%%%%%%%%%%%%%%%%%%%%%%%%%%%%%%%%%%%
\subsubsection{Ulam metric}
%%%%%%%%%%%%%%%%%%%%%%%%%%%%%%%%%%%%%%%%%%%%%%%%%%%%%%%%%%%%%%%%%%%%%%%%%%%%%%%%%

Let $L(\sigma,\tau)$ be the length of a \emph{longest common subsequence}. The \emph{Ulam distance} is
$$
 d_U(\sigma,\tau)=n-L(\sigma,\tau)$$
 % =n-\LIS(\sigma^{-1}\tau).$$
 
Equivalently, this is the minimum number of \emph{translocations}, each removing one symbol and reinserting it at another position~\cite{FarnoudSkachekMilenkovic2013}.

A {\emph{common subsequence}} is obtained by deleting symbols from each word without changing the order of those retained; its symbols need not occupy consecutive positions.

For example, $d_U(12345,14235)=1$, since $1235$ is a common subsequence and moving $4$ produces the second word.

%%%%%%%%%%%%%%%%%%%%%%%%%%%%%%%%%%%%%%%%%%%%%%%%%%%%%%%%%%%%%%%%%%%%%%%%%%%%%%%%%
\subsubsection{Kendall metric}
%%%%%%%%%%%%%%%%%%%%%%%%%%%%%%%%%%%%%%%%%%%%%%%%%%%%%%%%%%%%%%%%%%%%%%%%%%%%%%%%%

The \emph{Kendall distance} is the minimum number of swaps of adjacent positions needed to transform one permutation into the other, or equivalently
\[
 d_K(\sigma,\tau)=\inv(\sigma^{-1}\tau),\qquad
 \inv(\pi)=|\{(i,j):i<j,\ \pi(i)>\pi(j)\}|.
\]
Thus $d_K(1234,2134)=1$ and $d_K(12345,54321)=10$. This metric models local ordering errors~\cite{BuzagloEtzion2015}.

%%%%%%%%%%%%%%%%%%%%%%%%%%%%%%%%%%%%%%%%%%%%%%%%%%%%%%%%%%%%%%%%%%%%%%%%%%%%%%%%%
\subsubsection{Cayley metric}
%%%%%%%%%%%%%%%%%%%%%%%%%%%%%%%%%%%%%%%%%%%%%%%%%%%%%%%%%%%%%%%%%%%%%%%%%%%%%%%%%

The \emph{Cayley distance} is the minimum number of swaps of arbitrary pairs of positions needed to transform one permutation into the other. If $c(\pi)$ counts all cycles, including fixed points, then
\[
 d_C(\sigma,\tau)=n-c(\sigma^{-1}\tau).
\]
For example, $d_C(12345,21354)=2$, since the relative permutation is $(1\,2)(4\,5)$. The cycle formula follows by splitting each cycle into transpositions and observing that a transposition changes the number of cycles by one.

%%%%%%%%%%%%%%%%%%%%%%%%%%%%%%%%%%%%%%%%%%%%%%%%%%%%%%%%%%%%%%%%%%%%%%%%%%%%%%%%%
\subsection{Codes and their graph representations}
%%%%%%%%%%%%%%%%%%%%%%%%%%%%%%%%%%%%%%%%%%%%%%%%%%%%%%%%%%%%%%%%%%%%%%%%%%%%%%%%%
\begin{remark}
Under our composition convention, all four metrics are \emph{left-invariant}:
\[
 d_M(g\sigma,g\tau)=d_M(\sigma,\tau)=d_M(e,\sigma^{-1}\tau).
\]
The {Hamming and Cayley distances are} also \emph{right-invariant}; this is not true in general for Ulam or Kendall distance. The elementary operations give
\begin{align*}
 d_U(\sigma,\tau)&\le d_K(\sigma,\tau),\\d_C(\sigma,\tau)&\le d_K(\sigma,\tau),\\
 \tfrac12d_H(\sigma,\tau)\le d_C(\sigma,\tau)&\le d_H(\sigma,\tau)-1,\hbox{ for }\sigma\ne\tau\in S_n.
\end{align*}
For the last pair of inequalities, a nontrivial cycle of length $\ell$ contributes $\ell$ to Hamming distance and $\ell-1$ to Cayley distance.
\end{remark}

For $M\in\{H,U,K,C\}$ and an integer $d\ge1$, a \emph{permutation code} of minimum $d_M$-distance at least $d$ is a subset of $S_n$ whose distinct elements are pairwise at $d_M$-distance at least $d$. Let $A_M(n,d)$ denote its maximum possible size; explicitly,

\[
 A_M(n,d):=\max\{|\mathcal C|:\mathcal C\subseteq S_n,\ 
 d_M(\sigma,\tau)\ge d\text{ for distinct }\sigma,\tau\in\mathcal C\}.
\] 
If $d=2r+1$, the closed metric balls of radius $r$ centred at distinct codewords are disjoint: a point in two such balls would put their centres at distance at most $2r$, by the triangle inequality. This gives the sphere-packing bound~\eqref{eq:SP} used below (note that only this implication is needed; we will not assume the converse for an arbitrary metric on $S_n$).

A {\emph{finite simple undirected graph}} $G=(V,E)$ consists of a vertex set $V$ and a set of edges $E$; in particular, there are no loops or multiple edges. Vertices joined by an edge are \emph{adjacent}. An \emph{independent set} contains no adjacent pair, and its maximum possible size is the \emph{independence number} $\alpha(G)$.

For the coding problem of interest, let $\Gamma_M(n,d)$ have vertex set $S_n$, with distinct $\sigma$ and $\tau$ adjacent exactly when $d_M(\sigma,\tau)<d$. Thus each vertex represents a complete permutation, not an individual symbol, and an edge represents a pair forbidden by the minimum-distance requirement. It follows directly that
\[
 A_M(n,d)=\alpha\bigl(\Gamma_M(n,d)\bigr).
\]
Equivalently, a code is a {\emph{clique}}, that is, a set of pairwise adjacent vertices, in the \emph{complement} of $\Gamma_M(n,d)$, where adjacency is reversed for distinct vertices. We use the independent-set formulation throughout.

A \emph{path} is a sequence of distinct vertices with consecutive vertices adjacent; its length is the number of edges. The \emph{graph distance} between two vertices $u,v\in V(G)$, denoted $d_G(u,v)$, is the length of a shortest path from $u$ to $v$, with distance infinity when no such path exists. A graph is \emph{connected} when all pairs have finite distance; its \emph{diameter} is then the maximum of these distances.

For $M\in\{U,C,K\}$, let $G_M(n)$ join permutations at $d_M$-distance one. A sequence of elementary edits corresponds to a sequence of adjacent vertices, so the edit-distance definitions give
\[
 d_{G_M(n)}(\sigma,\tau)=d_M(\sigma,\tau).
\]
Consequently, a code of minimum distance $d$ has no path of length at most $d-1$ between distinct codewords in $G_M(n)$. For example, $123,213,231,321$ form a shortest path of length three in $G_K(3)$, since $\inv(321)=3$. The pair $\{123,321\}$ therefore has minimum Kendall distance three. The same distance-one construction is not suitable for Hamming distance: it gives an edgeless graph, whose graph distances do not recover $d_H$. {Section~\ref{sec:permutation-graphs} develops these graph representations and their parameters.}

%%%%%%%%%%%%%%%%%%%%%%%%%%%%%%%%%%%%%%%%%%%%%%%%%%%%%%%%%%%%%%%%%%%%%%%%%%%%%%%%%
\subsection{The Eigenvalue Method}\label{sec:eigenvalue-method}
%%%%%%%%%%%%%%%%%%%%%%%%%%%%%%%%%%%%%%%%%%%%%%%%%%%%%%%%%%%%%%%%%%%%%%%%%%%%%%%%%
{
The Eigenvalue Method~\cite{APR2025} combines polynomial spectral bounds for graph independence parameters with their corresponding optimization methods. We recall the inertia-type and ratio-type bounds below and describe the optimization of the polynomials appearing in these bounds. The inertia-type bound can be optimized by mixed-integer linear programming, while the ratio-type bound admits a linear programming formulation on partially walk-regular graphs. We refer to~\cite{BrouwerH2012} for the spectral graph theory background.
\par} Let $G=(V,E)$ be a finite simple undirected graph of {\emph{order}} $N=|V|$. Its \emph{adjacency matrix} $A$ has rows and columns indexed by $V$, with $A_{uv}=1$ when $u$ and $v$ are adjacent and $A_{uv}=0$ otherwise. Thus, $A$ is a real symmetric matrix with zero diagonal. {The \emph{degree} of a vertex is its number of neighbours.} A graph is \emph{$\delta$-regular} if every vertex has $\delta$ neighbours, in which case $A\one=\delta\one$, where $\one$ is the {\emph{all-one vector}}.

An \emph{automorphism} of $G$ is a permutation of its vertices that preserves adjacency. The graph is \emph{vertex-transitive} if, for any two vertices, some automorphism maps one to the other. For a group $\Gamma$ with identity element $e$, and an \emph{inverse-closed} set $S\subseteq\Gamma\setminus\{e\}$, meaning $S^{-1}=S$, the \emph{Cayley graph} $\Cay(\Gamma,S)$ joins $g$ and $h$ when $g^{-1}h\in S$. Its neighbours at $g$ are the elements $gs$ with $s\in S$. Left multiplication is an automorphism, so the graph is vertex-transitive and $|S|$-regular. It is connected precisely when $S$ generates $\Gamma$.

Since $A$ is real symmetric, it has an orthonormal basis of eigenvectors and real eigenvalues
\[
 \lambda_1\ge\lambda_2\ge\cdots\ge\lambda_N,
\]
listed with multiplicity. The \emph{spectrum} is this multiset of eigenvalues. We also use the distinct eigenvalues $\theta_0>\cdots>\theta_r$, with multiplicities $m_0,\ldots,m_r$, so that $\sum_{i=0}^r m_i=N$. In displayed spectra, $\theta_i^{m_i}$ indicates multiplicity, not exponentiation. For a connected $\delta$-regular graph, $\lambda_1=\theta_0=\delta$ is simple, with eigenvector $\one$, and all eigenvalues lie in $[-\delta,\delta]$. We call the remaining eigenvalues \emph{nonprincipal}; their eigenvectors lie in the subspace orthogonal to $\one$.

A \emph{walk} of length $j$ is a sequence of $j+1$ vertices with consecutive vertices adjacent; unlike a path, it may repeat vertices. Matrix multiplication shows that $(A^j)_{uv}$ counts walks of length $j$ from $u$ to $v$. In particular, $(A^j)_{vv}$ counts closed walks based at $v$, $(A^2)_{vv}$ is the degree of $v$, and $(A^3)_{vv}$ is twice the number of triangles containing $v$. We take $A^0=I$. {These observations explain the diagonal quantities in the inertia-type and ratio-type bounds below.}

A graph is \emph{bipartite} if its vertices can be partitioned into two independent sets. Every walk alternates between the two parts, so there are no closed walks of odd length. Its spectrum is symmetric about zero: changing the sign of an eigenvector on one part changes its eigenvalue from $\lambda$ to $-\lambda$. In particular, $(A^3)_{vv}=0$ for bipartite graphs, a simplification used in the cubic bounds below.

We say that a graph $G$ is \emph{$k$-partially walk-regular} when $(A^i)_{vv}$ is independent of $v$ for $0\le i\le k$. It is \emph{walk-regular} if this holds for every $i$. Vertex-transitivity implies walk-regularity. Indeed, an automorphism mapping $u$ to $v$ gives a bijection between closed walks of any given length based at these vertices. On a walk-regular graph, all diagonal entries of $p(A)$ are the same constant, for every polynomial $p$. Taking their sum and using the eigenvalues of $p(A)$ gives,
\begin{equation}\label{eq:diag}
 (p(A))_{vv}=\frac{1}{N}\tr p(A)
 =\frac1N\sum_{i=0}^r m_i p(\theta_i),\qquad N=|V|.
\end{equation}

\begin{definition}
The \emph{$k$-independence number} of a graph $G$, denoted $\alpha_k(G)$, is the largest cardinality of a vertex set whose distinct vertices are at graph distance greater than $k$. The \emph{graph power} $G^k$ joins distinct vertices at distance at most $k$; hence $\alpha_k(G)=\alpha(G^k)$ where $\alpha(G)$ denotes the usual \emph{independence number} of $G$.
\end{definition}

For the permutation base graphs $G_M(n)$, $M\in\{U,C,K\}$, a minimum-distance requirement $d\ge2$ corresponds to $k=d-1$. For Hamming distance we instead apply ordinary independence bounds, with $k=1$, directly to $\Gamma_H(n,d)$. The graph power $G^k$ should not be confused with the matrix power $A^k$: the latter counts walks and is generally not a zero--one adjacency matrix.

The spectra of $G$ and $G^k$ need not be related by an adjacency polynomial~\cite{AbiadCFNZ2022}. Applying a {polynomial spectral bound} to $\alpha_k(G)$ is therefore different from applying the ratio or inertia bound to $\alpha(G^k)$. {We compare the two approaches for Ulam distance in Table~\ref{tab:ulam}.}

%%%%%%%%%%%%%%%%%%%%%%%%%%%%%%%%%%%%%%%%%%%%%%%%%%%%%%%%%%%%%%%%%%%%%%%%%%%%%%%%%
%\subsubsection{The Eigenvalue Method:  inertia-type and ratio-type bounds}
%%%%%%%%%%%%%%%%%%%%%%%%%%%%%%%%%%%%%%%%%%%%%%%%%%%%%%%%%%%%%%%%%%%%%%%%%%%%%%%%%
{
We now recall the polynomial spectral bounds used in the Eigenvalue Method, together with the corresponding optimization methods; see~\cite{APR2025}.
\par} {We reserve \emph{ratio bound} for the degree-one case and use \emph{ratio-type bound} for the polynomial extensions.} 

Write $\mathbb R_k[x]$ for the real polynomials of degree at most $k$, where $k\ge1$. For $p(x)=\sum_{j=0}^k a_jx^j$, define the {\emph{matrix polynomial}} $p(A)=\sum_{j=0}^k a_jA^j$. It has eigenvalues $p(\lambda_i)$ with the same eigenvectors as $A$, although their ordering can change. The walk-counting interpretation gives the key implication
\[
 d_G(u,v)>k\quad\Longrightarrow\quad (p(A))_{uv}=0
 \qquad(u\ne v).
\]
Indeed, there is no walk of length at most $k$ between such vertices. If $\mathcal C$ is $k$-independent, the \emph{principal submatrix} obtained by restricting both the rows and the columns of $p(A)$ to $\mathcal C$ is therefore diagonal. On a walk-regular graph it is a scalar multiple of the identity. This converts the minimum-distance condition into a matrix constraint to which eigenvalue bounds apply~\cite{AbiadCF2019}. To allow for graphs without a common diagonal, put
\[
 W(p)=\max_v(p(A))_{vv},\qquad w(p)=\min_v(p(A))_{vv}.
\]
\begin{theorem}[Inertia-type bound~{\cite[Theorem 3.1]{AbiadCF2019}}]\label{th:inertia type bound}
For an $N$-vertex graph with eigenvalues $\lambda_1\ge\cdots\ge\lambda_N$, and $p\in\mathbb R_k[x]$,
\[
 \alpha_k(G)\le \min\bigl\{
 |\{i:p(\lambda_i)\ge w(p)\}|,
 |\{i:p(\lambda_i)\le W(p)\}|\bigr\}.
\]
\end{theorem}
%The mechanism is eigenvalue interlacing: eigenvalues of a principal submatrix lie between the corresponding ordered eigenvalues of the full symmetric matrix~\cite{H1995}. If a $k$-independent set has size $s$, its diagonal submatrix has $s$ eigenvalues in $[w(p),W(p)]$. Interlacing therefore forces at least $s$ eigenvalues of $p(A)$ to be at least $w(p)$ and at least $s$ to be at most $W(p)$, giving the two counts above.

For $k=1$ and $p(x)=x$ this is the Cvetkovi\'c \emph{inertia bound}~\cite{Cvetkovic1971}. On a walk-regular graph we may subtract the common diagonal of $p(A)$. If the resulting matrix has \emph{inertia} $(n_+,n_-,n_0)$, where these are the numbers of positive, negative and zero eigenvalues counted with multiplicity, respectively, the bound is $n_0+\min(n_+,n_-)$.

{
For the sharp inertia-type examples below, we use $p(x)=x^2-\delta$, whose diagonal entries vanish on a $\delta$-regular graph. We count its positive, negative and zero eigenvalues and exhibit codes attaining the resulting bounds. We therefore obtain the optimum in these cases without optimizing over all admissible polynomials~\cite{AbiadCFNZ2022}.
\par}

\begin{theorem}[Ratio-type bound~{\cite[Theorem 3.2]{AbiadCF2019}}]\label{th:ratio type bound}
Let $G$ be connected and $\delta$-regular, with $\lambda_1=\delta$. For $p\in\mathbb R_k[x]$, put $\lambda(p)=\min_{i\ge2}p(\lambda_i)$. If $p(\delta)>\lambda(p)$, then
\[
 \alpha_k(G)\le N\frac{W(p)-\lambda(p)}{p(\delta)-\lambda(p)}.
\]
\end{theorem}
%Here $\lambda(p)$ is the least eigenvalue of $p(A)$ on the subspace orthogonal to $\one$. Thus $y^Tp(A)y\ge\lambda(p)\|y\|^2$ whenever $y\perp\one$. Applying this inequality to the component of a code's characteristic vector orthogonal to $\one$, and using the diagonal submatrix described above, yields the ratio-type bound. Remark~\ref{rem:equality} gives this calculation and the equality condition.

For $k=1$ and $p(x)=x$, this gives Hoffman's \emph{ratio bound}; see, for example,~\cite{H1995}:
\begin{equation}\label{eq:ratio-one}
 \alpha(G)\le R_1(G):=-\frac{N\lambda_N}{\delta-\lambda_N}.
\end{equation}

{
The following formulas give the optimal ratio-type bounds in degrees two and three. For higher degrees, the polynomial optimization can be formulated as a linear program (LP)~\cite{Fiol2020}.
\par}

\begin{theorem}[Optimal quadratic ratio-type bound~{\cite[Corollary 3.3]{AbiadCF2019}}]\label{th: alpha_2}
Let $G$ be a connected $\delta$-regular graph on $N$ vertices, with distinct eigenvalues $\delta=\theta_0>\cdots>\theta_r$ and $r\ge2$. Let $\theta_i$ be the largest eigenvalue at most $-1$, and let $\theta_{i-1}$ be the next larger distinct eigenvalue. Then
\begin{equation}\label{eq:R2}
 \alpha_2(G)\le
 N\frac{\delta+\theta_i\theta_{i-1}}
 {(\delta-\theta_i)(\delta-\theta_{i-1})}.
\end{equation}
This is the optimum obtainable from Theorem~\ref{th:ratio type bound} with a polynomial of degree at most two.
\end{theorem}

\begin{theorem}[Optimal cubic ratio-type bound~{\cite[Theorem 11]{KN2023}}]\label{th: alpha3}
Let $G$ be a connected $\delta$-regular graph on $N$ vertices, with distinct eigenvalues $\delta=\theta_0>\cdots>\theta_r$, where $r\ge3$ and $\theta_r<-1$. Put $\Delta=\max_v(A^3)_{vv}$ and let $\theta_s$ be the smallest eigenvalue satisfying
\[
 \theta_s\ge-\frac{\delta^2+\delta\theta_r-\Delta}{\delta(\theta_r+1)}.
\]
Then the optimal degree-at-most-three ratio-type bound is
\begin{equation}\label{eq:R3}
 \alpha_3(G)\le N\frac{\Delta-\delta(\theta_s+\theta_{s+1}+\theta_r)
 -\theta_s\theta_{s+1}\theta_r}
 {(\delta-\theta_s)(\delta-\theta_{s+1})(\delta-\theta_r)}.
\end{equation}
\end{theorem}
The optimizing polynomials in~\eqref{eq:R2} and~\eqref{eq:R3} have the displayed two and three roots, respectively. Degenerate small cases can be evaluated directly using Theorem~\ref{th:ratio type bound}.

For general $k$, on a connected walk-regular graph (which is necessarily regular), the normalized minor-polynomial LP of Fiol~\cite{Fiol2020} is equivalently
\begin{equation}\label{eq: lp ratio}
 \operatorname{minimize}_{f\in\mathbb R_k[x]}
 \sum_{j=0}^r m_jf(\theta_j),\qquad
 f(\delta)=1,\quad f(\theta_j)\ge0\ (j\ge1).
\end{equation}
To see the normalization, set $f(x)=(p(x)-\lambda(p))/(p(\delta)-\lambda(p))$ for an admissible $p$. Then $f(\delta)=1$ and the other eigenvalues of $f(A)$ are nonnegative. Such a symmetric matrix is called \emph{positive semidefinite} (PSD). By~\eqref{eq:diag}, its trace is exactly $N(W(p)-\lambda(p))/(p(\delta)-\lambda(p))$, the corresponding ratio-type bound. The above is a linear program in the coefficients of $f$. It is equivalent to the divided-difference formulation in~\cite{Fiol2020}. Its optimum is the best ratio-type bound (we derive integer code bounds by taking floors).

{
\begin{remark}[Equality in the ratio-type bound]\label{rem:equality}
On a connected walk-regular graph, put $a=(p(A))_{vv}$ and $m=\lambda(p)$. For a $k$-independent set $\mathcal C$ of size $s$, its characteristic vector $z$ satisfies $z^Tp(A)z=sa$. Write $z=(s/N)\one+y$, where $y\perp\one$. Then
\[
 sa\ge\frac{s^2}{N}p(\delta)+m\left(s-\frac{s^2}{N}\right).
\]
For nonempty $\mathcal C$, equality in the real-valued ratio-type bound holds exactly when $p(A)y=my$. Since code sizes are integers, the floor of a ratio-type bound may be attained even when the displayed inequality is strict.
\end{remark}
\par}

%%%%%%%%%%%%%%%%%%%%%%%%%%%%%%%%%%%%%%%%%%%%%%%%%%%%%%%%%%%%%%%%%%%%%%%%%%%%%%%%%
\section{Known bounds for permutation codes}\label{sec:known-bounds}
%%%%%%%%%%%%%%%%%%%%%%%%%%%%%%%%%%%%%%%%%%%%%%%%%%%%%%%%%%%%%%%%%%%%%%%%%%%%%%%%%

The upper bounds we will be using for comparison with the spectral bounds come from packing and anticode arguments, injective projections giving Singleton-type bounds, and algebraic optimization. The refinements available depend on the specific permutation metric. For the Hamming distance, linear programming in the conjugacy scheme of $S_n$~\cite{Tarnanen1999,Bogaerts2010}, semidefinite programming~\cite{BogaertsDukes2014}, and combinatorial arguments~\cite{BierbrauerMetsch2013} provide upper bounds. For the Ulam distance, the subsequence description gives a Singleton bound, while linear and integer programming give further upper bounds, and computer searches determine exact small values~\cite{FarnoudSkachekMilenkovic2013,Riet}. For the Cayley distance, the dependence on cycle type permits the conjugacy-class LP described below. For Kendall distance, anticode refinements~\cite{BuzagloEtzion2015}, coherent configurations~\cite{LimHagiwara2012}, and representation-theoretic integer programming~\cite{Abdollahi2023PermutationCodes} give useful comparisons.

A common starting point is the \emph{code--anticode bound}. An {\emph{anticode}} of diameter at most $d-1$ is a nonempty subset $\mathcal A\subseteq S_n$ whose pairwise $d_M$-distances are at most $d-1$. For each of the four metrics,
\begin{equation}\label{eq:code-anticode}
 A_M(n,d)\le\left\lfloor\frac{n!}{|\mathcal A|}\right\rfloor.
\end{equation}
Indeed, left invariance implies that a code $\mathcal C$ of minimum distance at least $d$ meets each translate $g\mathcal A$ in at most one element. Summing $|\mathcal C\cap g\mathcal A|$ over $g\in S_n$ gives $|\mathcal C|\,|\mathcal A|\le n!$. This argument underlies the Kendall bounds from balls and double balls in~\cite{BuzagloEtzion2015}.

Throughout, we use the term \emph{sphere-packing bound} for the disjoint-ball bound~\eqref{eq:SP}, with closed balls of radius $t=\lfloor(d-1)/2\rfloor$ in the metric under consideration. We give their sizes $B_M(n,t)$ in~\eqref{eq:ballsH}--\eqref{eq:ballsU}. Such a ball has diameter at most $d-1$, so sphere packing is a special case of~\eqref{eq:code-anticode}. The following subsections state the relevant Singleton bounds and the more specialized bounds used in our tables. Tables~\ref{tab:hamming}, \ref{tab:cayley}, and~\ref{tab:ulam} explicitly include the corresponding Singleton comparisons; the Kendall comparisons use the stronger small-distance bounds described below.

%%%%%%%%%%%%%%%%%%%%%%%%%%%%%%%%%%%%%%%%%%%%%%%%%%%%%%%%%%%%%%%%%%%%%%%%%%%%%%%%%
\subsection{Hamming distance}
%%%%%%%%%%%%%%%%%%%%%%%%%%%%%%%%%%%%%%%%%%%%%%%%%%%%%%%%%%%%%%%%%%%%%%%%%%%%%%%%%

For $1\le d\le n$, the classical permutation \emph{Singleton bound} is~\cite[Theorem 1]{Bogaerts2010}
\begin{equation}\label{eq:H-singleton}
 A_H(n,d)\le\frac{n!}{(d-1)!}.
\end{equation}
To see this, retain any fixed $n-d+1$ positions of each codeword. Two codewords with the same retained entries would have Hamming distance at most $d-1$, so this projection is injective. It has $n!/(d-1)!$ possible values.

Stronger bounds include the linear programming bounds of Tarnanen and Bogaerts~\cite{Tarnanen1999,Bogaerts2010}, semidefinite programming bounds~\cite{BogaertsDukes2014}, and combinatorial bounds~\cite{BierbrauerMetsch2013}. Exact small values can also be obtained by clique methods~\cite{Barta2022}. For $(n,d)=(7,5)$, Montemanni, Barta and Smith~\cite{Montemanni2014} give the upper bound $122$, improving the semidefinite programming bound $134$.

%%%%%%%%%%%%%%%%%%%%%%%%%%%%%%%%%%%%%%%%%%%%%%%%%%%%%%%%%%%%%%%%%%%%%%%%%%%%%%%%%
\subsection{Ulam distance}
%%%%%%%%%%%%%%%%%%%%%%%%%%%%%%%%%%%%%%%%%%%%%%%%%%%%%%%%%%%%%%%%%%%%%%%%%%%%%%%%%

For $1\le d\le n$, the \emph{Singleton bound} is~\cite{FarnoudSkachekMilenkovic2013}
\begin{equation}\label{eq:U-singleton}
 A_U(n,d)\le(n-d+1)!.
\end{equation}
Here one retains a fixed set of $n-d+1$ symbols in their original order. Equal projected words would give a common subsequence of length $n-d+1$, contrary to minimum Ulam distance $d$. Thus the projection is injective and has $(n-d+1)!$ possible values.

Farnoud, Skachek and Milenkovic~\cite{FarnoudSkachekMilenkovic2013} develop the coding framework, while G\"olo\u{g}lu, Lember, Riet and Skachek~\cite{Riet} give linear/integer programming bounds and small exact values. Recent constructions are considered in~\cite{10.1007/978-3-032-17156-6_4}. We use the exact values in~\cite[Table I]{Riet} for comparison.

%%%%%%%%%%%%%%%%%%%%%%%%%%%%%%%%%%%%%%%%%%%%%%%%%%%%%%%%%%%%%%%%%%%%%%%%%%%%%%%%%
\subsection{Cayley distance}
%%%%%%%%%%%%%%%%%%%%%%%%%%%%%%%%%%%%%%%%%%%%%%%%%%%%%%%%%%%%%%%%%%%%%%%%%%%%%%%%%

A Singleton-type projection argument gives
\begin{equation}\label{eq:C-singleton}
 A_C(n,d)\le\frac{n!}{d!},\qquad 1\le d\le n.
\end{equation}
Retain any fixed $n-d$ positions. Two permutations with the same retained entries differ only on the remaining $d$ positions, so their relative permutation has Cayley length at most $d-1$. The projection is therefore injective on a code of minimum distance at least $d$, and it has $n!/d!$ possible values.

The sphere-packing bound~\eqref{eq:SP} and the conjugacy-class linear programming bound~\eqref{eq:cayley-character-lp} give further upper bounds. The latter uses the symmetric-group character framework of Tarnanen~\cite{Tarnanen1999}; the formulation allowing arbitrary admissible conjugacy classes is stated in~\cite[Theorem 6]{Bogaerts2010}. We specify its objective and constraints in~\eqref{eq:cayley-character-lp}, selecting the admissible classes by Cayley distance, and compare its values with the spectral, sphere-packing, and Singleton bounds in Table~\ref{tab:cayley}.

For minimum Cayley distance three, a radius-one ball consists of its centre and the $\binom n2$ transpositions. Thus $B_C(n,1)=1+\binom n2$, and~\eqref{eq:SP} becomes
\begin{equation}\label{eq:CSP3}
 A_C(n,3)\le\left\lfloor\frac{n!}{1+\binom n2}\right\rfloor.
\end{equation}
The explicit spectrum of the transposition graph~\cite{KonstantinovaKravchuk2022} makes the Eigenvalue Method particularly accessible. Finally, note that Nguyen~\cite{Nguyen2024GilbertVarshamov} derived improved lower bounds.

%%%%%%%%%%%%%%%%%%%%%%%%%%%%%%%%%%%%%%%%%%%%%%%%%%%%%%%%%%%%%%%%%%%%%%%%%%%%%%%%%
\subsection[Kendall-tau distance]{Kendall-$\tau$ distance}
%%%%%%%%%%%%%%%%%%%%%%%%%%%%%%%%%%%%%%%%%%%%%%%%%%%%%%%%%%%%%%%%%%%%%%%%%%%%%%%%%

For $n-1<d\le\binom n2$, Barg and Mazumdar's Singleton bound~\cite[Theorem 3.2]{BargMazumdar2010} gives
\begin{equation}\label{eq:K-singleton}
 A_K(n,d)\le
 \left\lfloor\frac32+\sqrt{n(n-1)-2d+\frac14}\right\rfloor !.
\end{equation}
It follows by retaining the relative order of a fixed set of symbols. This projection bound improves the trivial bound $n!$ only when $d>n-1$, so it gives no improvement for the rows with $n\ge5$ in Table~\ref{tab:kendall}.

For small even distances, the double-ball bound is a useful refinement of sphere packing~\cite{BuzagloEtzion2015}. Write $\mathcal B_K(\sigma,t)$ for the closed Kendall ball of radius $t$ centred at $\sigma$. If $d_K(e,s)=1$, then $\mathcal B_K(e,t)\cup\mathcal B_K(s,t)$ has diameter at most $2t+1$. Applying~\eqref{eq:code-anticode} to this union gives
\begin{equation}\label{eq:K-double-ball}
 A_K(n,2t+2)\le
 \left\lfloor\frac{n!}{|\mathcal B_K(e,t)\cup\mathcal B_K(s,t)|}\right\rfloor.
\end{equation}
For $t=1$, each ball has size $n$ and their intersection consists of $e$ and $s$, since the adjacent-transposition graph is bipartite. Hence
\begin{equation}\label{eq:K-double-ball-four}
 A_K(n,4)\le\left\lfloor\frac{n!}{2(n-1)}\right\rfloor,
\end{equation}

which gives $420$ and $2880$ at $n=7$ and $n=8$, respectively. This is the bound of~\cite[Corollary~6]{BuzagloEtzion2015}. Moreover,~\cite[Theorem~8]{BuzagloEtzion2015} shows that, for $n\ge4$, a double ball of radius one is a largest Kendall anticode of diameter three. Thus~\eqref{eq:K-double-ball-four} is the strongest upper bound obtained from~\eqref{eq:code-anticode} using a diameter-three anticode.

Further upper bounds use coherent-configuration methods~\cite{LimHagiwara2012} and representation-theoretic integer programming~\cite{Abdollahi2023PermutationCodes}. The comparison incorporates the updated table of Parvaresh et al.~\cite{ParvareshEtAl2025} and the exact small values reported in~\cite{Vijayakumaran2016}. In particular, the bounds $A_K(6,3)\le116$ and $A_K(7,3)\le716$ of~\cite{Abdollahi2023PermutationCodes} are {stronger than sphere packing}.

%%%%%%%%%%%%%%%%%%%%%%%%%%%%%%%%%%%%%%%%%%%%%%%%%%%%%%%%%%%%%%%%%%%%%%%%%%%%%%%%%
\section{Graphs associated to permutation codes}\label{sec:permutation-graphs}
%%%%%%%%%%%%%%%%%%%%%%%%%%%%%%%%%%%%%%%%%%%%%%%%%%%%%%%%%%%%%%%%%%%%%%%%%%%%%%%%%
We distinguish a graph defined by a minimum-distance threshold from a distance-one base graph. For a finite integer-valued metric space $(X,d_X)$, we can define a graph with 
\[
 V(G)=X,\qquad x\sim_G y\ \Longleftrightarrow\ d_X(x,y)=1.
\]
When $d_X$ agrees with the shortest-path distance in $G$, we identify its $k$-th power $G^k$ with the graph joining distinct $x,y$ at metric distance at most $k$. In this case, a minimum-distance-$(k+1)$ code is a $k$-independent set in $G$~\cite{APR2025}. For permutations we write $G_M(n)$ for the base graph and $\Gamma_M(n,d)$ for the forbidden-distance graph. The identity between the latter and $G_M(n)^{d-1}$ holds for $M\in\{K,C,U\}$; the Hamming exception is explained in Remark~\ref{rem:hamming}.

%%%%%%%%%%%%%%%%%%%%%%%%%%%%%%%%%%%%%%%%%%%%%%%%%%%%%%%%%%%%%%%%%%%%%%%%%%%%%%%%%

%%%%%%%%%%%%%%%%%%%%%%%%%%%%%%%%%%%%%%%%%%%%%%%%%%%%%%%%%%%%%%%%%%%%%%%%%%%%%%%%%

\subsection[Forbidden-distance graphs (power graphs) and independence bounds]{The forbidden-distance graph  (power graph) $\Gamma_M(n,d)$ and bounds on $\alpha(\Gamma_M(n,d))$}
%%%%%%%%%%%%%%%%%%%%%%%%%%%%%%%%%%%%%%%%%%%%%%%%%%%%%%%%%%%%%%%%%%%%%%%%%%%%%%%%%
For $M\in\{H,U,K,C\}$, the graph used to encode the minimum-distance requirement is
\[
 \Gamma_M(n,d):\quad V=S_n,\qquad
 \sigma\sim\tau\ \Longleftrightarrow\ 0<d_M(\sigma,\tau)<d.
\]
Then $A_M(n,d)=\alpha(\Gamma_M(n,d))$. Barta and Montemanni~\cite{Barta2022} study permutation Hamming graphs joining pairs with $d_H\ge d$, while G\"olo\u{g}lu et al.~\cite{Riet} use the analogous graph with $d_U\ge d$ for Ulam-code searches. Their underlying unweighted graphs are therefore the complements of $\Gamma_H(n,d)$ and $\Gamma_U(n,d)$, respectively. For $M\in\{K,C,U\}$, in Proposition~\ref{prop:basegraphs} below we identify $\Gamma_M(n,d)$ with $G_M(n)^{d-1}$.

\begin{lemma}\label{lem:regular}
Every $\Gamma_M(n,d)$ is a Cayley graph, hence vertex-transitive and walk-regular. If $B_M(n,t):=|\{\pi:d_M(e,\pi)\le t\}|$, its degree is $B_M(n,d-1)-1$.
\end{lemma}
\begin{proof}
Left invariance makes $\pi\mapsto\sigma^{-1}\pi$ an isometry sending $\sigma$ to $e$. Thus $\Gamma_M(n,d)$ is the Cayley graph with connecting set $\{\pi:0<d_M(e,\pi)<d\}$. For any $\pi$, left invariance and symmetry give $d_M(e,\pi^{-1})=d_M(\pi,e)=d_M(e,\pi)$, so this set is inverse-closed. Its size is $B_M(n,d-1)-1$, since the ball contains its centre. Left translations act transitively on the vertices; they also identify the closed walks based at any two vertices, proving walk-regularity.

For $d\ge2$, the Ulam, Cayley and Kendall connecting sets contain all adjacent transpositions, so these graphs are connected. For Hamming distance the graph is connected for $d\ge3$, since an adjacent transposition has Hamming distance two. At $d=2$, $\Gamma_H(n,2)$ is edgeless.
\end{proof}

\allowdisplaybreaks

For reference, the ball sizes are
\begin{align}
 B_H(n,t)&=\sum_{j=0}^t\binom nj D_j,&D_j&=j!\sum_{i=0}^j\frac{(-1)^i}{i!},\label{eq:ballsH}\\
 B_C(n,t)&=\sum_{j=0}^t[z^j]\prod_{i=1}^{n-1}(1+iz),\label{eq:ballsC}\\
 B_K(n,t)&=\sum_{j=0}^t[z^j]\prod_{i=1}^{n}(1+z+\cdots+z^{i-1}),\label{eq:ballsK}\\
 B_U(n,t)&=\sum_{\substack{\lambda\vdash n\\\lambda_1\ge n-t}}(f^\lambda)^2.\label{eq:ballsU}
\end{align}
Here $D_j$ is the number of derangements of $j$ symbols, $[z^j]$ denotes coefficient extraction, and $f^\lambda$ is the number of standard Young tableaux of shape $\lambda$. The notation $\lambda\vdash n$ means that $\lambda=(\lambda_1,\ldots,\lambda_r)$ is a {\emph{partition}} of $n$, with $\lambda_1\ge\cdots\ge\lambda_r>0$ and $\sum_{i=1}^r\lambda_i=n$. Its {\emph{Young diagram}} has $\lambda_i$ boxes in row $i$; in particular, $\lambda_1$ is the largest part and the length of the first row. A {\emph{standard Young tableau}} is a filling of these boxes with $1,\ldots,n$, each used once, increasing along rows and down columns; see~\cite[Chapters 2--3]{Sagan2001}. The first three identities count by support, cycles and inversions; the last follows from the Robinson--Schensted correspondence~\cite{Sagan2001,Riet}. Under this correspondence, $\lambda_1$ is the length of a longest increasing subsequence of the permutation, so $d_U(e,\pi)\le t$ is equivalent to $\lambda_1\ge n-t$. Disjoint balls give the \emph{sphere-packing bound}
\begin{equation}\label{eq:SP}
 A_M(n,d)\le\left\lfloor\frac{n!}{B_M(n,\lfloor(d-1)/2\rfloor)}\right\rfloor.
\end{equation}

%%%%%%%%%%%%%%%%%%%%%%%%%%%%%%%%%%%%%%%%%%%%%%%%%%%%%%%%%%%%%%%%%%%%%%%%%%%%%%%%%
\subsection[Base graph and higher independence bounds]{The base graph $G_M(n)$ and bounds on $\alpha_k(G_M(n))$}\label{sec:basegraphs}
%%%%%%%%%%%%%%%%%%%%%%%%%%%%%%%%%%%%%%%%%%%%%%%%%%%%%%%%%%%%%%%%%%%%%%%%%%%%%%%%%

Here $G_M(n)$ denotes the graph on $S_n$ joining permutations at $d_M$-distance one. For an inverse-closed set $S\subseteq S_n\setminus\{e\}$, the notation $\Cay(S_n,S)$ means the graph with vertex set $S_n$ and edges $\{\sigma,\sigma s\}$, for $\sigma\in S_n$ and $s\in S$. Equivalently, $\sigma\sim\tau$ when $\sigma^{-1}\tau\in S$. With our composition convention, right multiplication by a generator acts on positions. The three base graphs are
\begin{align*}
 G_K(n)&=\Cay(S_n,\{(i\ i+1):1\le i<n\}),\\
 G_C(n)&=\Cay(S_n,\{(i\ j):1\le i<j\le n\}),\\
 G_U(n)&=\Cay(S_n,T\cup T^{-1}),\quad
 T=\{(a\ a+1\ \cdots\ b):1\le a<b\le n\}.
\end{align*}
The first is the adjacent-transposition graph associated with Kendall distance, used for example by Buzaglo and Etzion~\cite{BuzagloEtzion2015}; the second is the all-transposition graph studied spectrally by Konstantinova and Kravchuk~\cite{KonstantinovaKravchuk2022}. The third uses the translocations defining Ulam distance~\cite{FarnoudSkachekMilenkovic2013}. The cycle $(a\ a+1\ \cdots\ b)$ moves the entry in position $a$ to position $b$, shifting the intervening entries left; its inverse moves the entry in position $b$ to position $a$. The union $T\cup T^{-1}$ is a set: adjacent transpositions are self-inverse and are counted only once.

\begin{proposition}\label{prop:basegraphs}
For $M\in\{K,C,U\}$, the graph distance in $G_M(n)$ equals $d_M$, and
\[
 \Gamma_M(n,d)=G_M(n)^{d-1},\qquad
 A_M(n,d)=\alpha_{d-1}(G_M(n)).
\]
These graphs are connected, regular and walk-regular, with parameters
\begin{center}
\begin{tabular}{lccc}
\toprule
Metric & Vertices & Degree & Diameter\\
\midrule
Kendall & $n!$ & $n-1$ & $\binom n2$\\
Cayley & $n!$ & $\binom n2$ & $n-1$\\
Ulam & $n!$ & $(n-1)^2$ & $n-1$\\
\bottomrule
\end{tabular}
\end{center}
\end{proposition}
\begin{proof}
 For Ulam distance, a common subsequence of length $L$ can be left fixed while the other $n-L$ symbols are moved, and one translocation can increase the LCS by at most one. Thus its graph distance is $n-L$. The Ulam degree is $(n-1)+2\sum_{\ell=3}^n(n-\ell+1)=(n-1)^2$. The other degrees and the diameters follow from the inversion and cycle formulas. All connecting sets contain the adjacent transpositions, so the graphs are connected. The remaining properties follow from Lemma~\ref{lem:regular}.
\end{proof}

\begin{remark}\label{rem:hamming}
Let $\mathcal H_n=H(n,n)$ be the Hamming graph on $[n]^n$, whose edges join words differing in exactly one coordinate. Its shortest-path distance is the Hamming distance. Identifying permutations with their one-line words gives $S_n\subseteq[n]^n$, and hence, for $d\ge2$,
\begin{equation}\label{eq:hamming-ambient-power}
 \Gamma_H(n,d)=\bigl(\mathcal H_n^{\,d-1}\bigr)[S_n].
\end{equation}
Here $X[S_n]$ denotes the subgraph of $X$ induced by the permutation vertices. Thus the graph power is formed before restricting to $S_n$. Restricting first gives the edgeless graph $\mathcal H_n[S_n]$, because distinct permutations cannot differ in exactly one coordinate. Its shortest-path distances do not recover $d_H$, and its powers differ from~\eqref{eq:hamming-ambient-power} for $d\ge3$.

In this case, $A_H(n,d)=\alpha(\Gamma_H(n,d))$. A valid (but generally weaker) use of the Cayley base graph follows from $d_H\le2d_C$:
\[
 A_H(n,d)\le A_C(n,\lceil d/2\rceil)
 =\alpha_{\lceil d/2\rceil-1}(G_C(n)).
\]
In particular, replacing ``distance one'' by ``Hamming distance two'' gives the Cayley metric, not Hamming distance divided by two. For example, $1234$ and $2143$ have Hamming distance four and Cayley distance two. The path $1234,2134,2143$ shows that they are adjacent in $G_C(4)^2$, whereas they are not adjacent in $\Gamma_H(4,4)$.
\end{remark}

%\paragraph{Spectral input.}
For each partition $\lambda\vdash n$, let $\rho^\lambda$ be the corresponding \emph{irreducible complex representation} of $S_n$, of dimension $f^\lambda$~\cite{jamesrepresentation}. For an inverse-closed connecting set $S$, the \emph{spectrum of $\Cay(S_n,S)$} is the union of the spectra of
\begin{equation}\label{eq:blocks}
 B_\lambda=\sum_{s\in S}\rho^\lambda(s),\qquad\lambda\vdash n,
\end{equation}
with every block eigenvalue repeated $f^\lambda$ times, by the decomposition of the regular representation~\cite{Sagan2001}. We use rational {\emph{Young seminormal matrices}}~\cite{Sagan2001}. To specify the matrices used in the computations, index a basis $\{v_t\}$ by standard Young tableaux $t$ of shape $\lambda$, let $c_t(j)$ be the column index minus the row index of the box containing $j$, and set $a_i(t)=c_t(i+1)-c_t(i)$. For $s_i=(i\ i+1)$, the action is
\[
 \rho^\lambda(s_i)v_t=\frac{1}{a_i(t)}v_t+
 \left(1+\frac{1}{a_i(t)}\right)v_{s_it},
\]
where the second term is omitted if $s_it$ is not standard. Here $s_it$ is obtained by interchanging the entries $i$ and $i+1$. We form products of these matrices to represent the other generators in~\eqref{eq:blocks} and compute the characteristic polynomials of the resulting blocks over $\mathbb Q$.

For a partition $\mu\vdash n$, let $\mathcal K_\mu$ be the {\emph{conjugacy class}} consisting of all permutations whose cycle lengths, including fixed points, form $\mu$, and let $|\mathcal K_\mu|$ denote its size. The {\emph{irreducible character}} associated with $\rho^\lambda$ is $\chi^\lambda(g)=\tr\rho^\lambda(g)$ for $g\in S_n$. Characters are constant on conjugacy classes, so we write $\chi^\lambda(\mu)$ for this value at any $g\in\mathcal K_\mu$; in particular, $\chi^\lambda(e)=f^\lambda$. See~\cite[Chapters 1--2]{Sagan2001} for these representation-theoretic definitions.

If the connecting set $S$ is a union of conjugacy classes, Schur's lemma makes each block $B_\lambda$ a scalar matrix. Taking its trace gives the associated eigenvalue~\cite[Proposition 6.3.1]{BrouwerH2012}:
\begin{equation}\label{eq:chars}
 \theta_\lambda=\frac{1}{f^\lambda}
 \sum_{\substack{\mu\vdash n\\\mathcal K_\mu\subseteq S}}
 |\mathcal K_\mu|\chi^\lambda(\mu),
 \qquad\lambda\vdash n.
\end{equation}
Thus the sum runs over the cycle types whose entire conjugacy classes lie in $S$. Each partition $\lambda$ contributes $(f^\lambda)^2$ copies of $\theta_\lambda$ to the spectrum; these multiplicities are added when different partitions give the same eigenvalue. For an explicit statement of this character formula and its multiplicity rule, see~\cite[Theorem 3]{KonstantinovaKravchuk2022}.

This applies to $\Gamma_H(n,d)$ and $\Gamma_C(n,d)$: both connecting sets are unions of conjugacy classes, since support size and cycle count depend only on cycle type. In particular,
\begin{equation}\label{eq:contents}
 \Spec G_C(n)=\left\{\left(\sum_{(i,j)\in\lambda}(j-i)\right)^{(f^\lambda)^2}:\lambda\vdash n\right\},
\end{equation}
with coincident eigenvalues combined~\cite{KonstantinovaKravchuk2022}. The Kendall and Ulam generating sets are not conjugacy-invariant, so one must use the full blocks in~\eqref{eq:blocks}.% In particular, a weighted random-to-random shuffle matrix is not the unweighted adjacency matrix of $G_U(n)$.

%%%%%%%%%%%%%%%%%%%%%%%%%%%%%%%%%%%%%%%%%%%%%%%%%%%%%%%%%%%%%%%%%%%%%%%%%%%%%%%%%
\section{The spectral upper bounds}\label{sec:new-bounds}
%%%%%%%%%%%%%%%%%%%%%%%%%%%%%%%%%%%%%%%%%%%%%%%%%%%%%%%%%%%%%%%%%%%%%%%%%%%%%%%%%

We write $R_k(G)$ for the optimal ratio-type bound obtained from polynomials of degree at most $k$. In particular, $R_1(G)$ is the ratio bound~\eqref{eq:ratio-one}. We evaluate~\eqref{eq:ratio-one}, \eqref{eq:R2}, and~\eqref{eq:R3} using the spectra in~\eqref{eq:blocks}--\eqref{eq:contents}, and take the floor to obtain integer upper bounds on $A_M(n,d)$.

In the tables, bold entries give known exact values. A superscript $\dagger$ indicates that a bound is attained only after taking the floor, whereas a bold bound without $\dagger$ is attained before rounding. A superscript $\star$ indicates an improvement over the cited best known upper bound. The Hamming and Kendall tables distinguish exact values from previously known upper bounds in the last column. We also prove that the inertia-type bound is sharp in Propositions~\ref{prop:C43} and~\ref{prop:K33}.

%%%%%%%%%%%%%%%%%%%%%%%%%%%%%%%%%%%%%%%%
\subsection{Theoretical comparison criteria}\label{sec:low-degree-comparisons}
%%%%%%%%%%%%%%%%%%%%%%%%%%%%%%%%%%%%%%%%
For $M\in\{U,C,K\}$, the quadratic ratio-type bound and the cubic ratio-type bound concern minimum distances three and four, respectively. The Hamming ratio-type bounds below instead use degree one on the forbidden-distance graph. We first give an exact comparison for the quadratic ratio-type bound. On an $N$-vertex $\delta$-regular base graph, a radius-one ball has $\delta+1$ vertices, so the real-valued sphere-packing benchmark for minimum distance three is $S_2=N/(\delta+1)$. We write $\Spec(G)$ for the multiset of eigenvalues of the adjacency matrix of $G$.

\begin{proposition}[Quadratic ratio-type bound: improvement criterion]\label{prop:quadratic-comparison}
Under the hypotheses of Theorem~\ref{th: alpha_2}, put $u=\theta_i\le-1$ and $v=\theta_{i-1}>-1$. Then
\begin{equation}\label{eq:quadratic-gap}
 g_2:=\frac{N}{\delta+1}-R_2(G)
 =-\frac{N\delta(u+1)(v+1)}
 {(\delta+1)(\delta-u)(\delta-v)}\ge0.
\end{equation}
Thus the optimal quadratic ratio-type bound is no weaker than sphere packing before rounding, and the inequality is strict if and only if $-1\notin\Spec(G)$. For comparison with a known bound, let $B$ be a fixed integer such that $\alpha_2(G)\le B$. Then
\begin{equation}\label{eq:quadratic-integer-comparison}
 \lfloor R_2(G)\rfloor<B
 \quad\Longleftrightarrow\quad
 g_2>\frac{N}{\delta+1}-B.
\end{equation}
In particular, strict improvement over the integer sphere-packing bound occurs exactly when $g_2>S_2-\lfloor S_2\rfloor$.
\end{proposition}
\begin{proof}
Subtracting~\eqref{eq:R2} from $N/(\delta+1)$ gives~\eqref{eq:quadratic-gap}. The denominator is positive, $u+1\le0$, and $v+1>0$, so $g_2=0$ exactly when $u=-1$. Since $u$ is the largest eigenvalue at most $-1$, this is equivalent to $-1\in\Spec(G)$. Finally, $\lfloor R_2(G)\rfloor<B$ is equivalent to $R_2(G)<B$ because $B$ is an integer, proving~\eqref{eq:quadratic-integer-comparison}.
\end{proof}

The previous results shows that, for a fixed length, a comparison of the spectral bound with a specified published coding-theoretic bound is decidable by exact finite calculations. %The adjacency matrix is an integer matrix; alternatively, the rational blocks in~\eqref{eq:blocks} give characteristic polynomials over $\mathbb Q$. Their real roots can be isolated and compared exactly, including equality cases, and the resulting quadratic ratio-type bound or cubic ratio-type bound can be compared with the chosen integer benchmark. For the cubic ratio-type bound in the Kendall metric, Corollary~\ref{cor:kendall-cubic-comparison} reduces this test to a spectral-gap inequality; against the double-ball bound, only nonsingularity needs to be checked.  

%%%%%%%%%%%%%%%%%%%%%%%%%%%%%%%%%%%%%%%%%%%%%%%%%%%%%%%%%%%%%%%%%%%%%%%%%%%%%%%%%
\subsection{Hamming metric}
%%%%%%%%%%%%%%%%%%%%%%%%%%%%%%%%%%%%%%%%%%%%%%%%%%%%%%%%%%%%%%%%%%%%%%%%%%%%%%%%%
We work with $\Gamma_H(n,d)$, whose vertices are permutations and whose edges join pairs with $0<d_H<d$. It is the complement of the underlying permutation Hamming graph used in the clique approach of Barta and Montemanni~\cite{Barta2022}. We obtain Table~\ref{tab:hamming} by applying the ratio bound directly to $\Gamma_H(n,d)$. We compare these bounds with known exact values and previously known upper bounds in the last two columns.
\begin{table}[htbp]
\centering\small
\caption{Upper bounds on $A_H(n,d)$ from the ratio bound~\eqref{eq:ratio-one} and the Singleton bound~\eqref{eq:H-singleton}, together with known exact values and previously known upper bounds.}\label{tab:hamming}

\setlength{\tabcolsep}{4pt}
\begin{tabular}{@{}lrrrl@{}}
\toprule
$(n,d)$ & \makecell{Ratio bound\\\eqref{eq:ratio-one}} & \makecell{Singleton\\\eqref{eq:H-singleton}} & \makecell{Known value or\\upper bound} & Reference\\
\midrule
$(4,4)$ & 5 & \textbf{4} & \textbf{4} & Exact; proof below\\
$(5,4)$ & \textbf{20} & \textbf{20} & \textbf{20} & Exact; Proposition~\ref{prop:hsharp}\\
$(6,4)$ & \textbf{120} & \textbf{120} & \textbf{120} & Exact~\cite{Barta2022}\\
$(6,5)$ & 68 & 30 & \textbf{18} & Exact~\cite{Barta2022}\\
$(7,4)$ & 543 & 840 & 535 & Upper bound (SDP)~\cite{BogaertsDukes2014}\\
$(7,5)$ & 610 & 210 & 122 & Upper bound (combinatorial)~\cite{Montemanni2014}\\
\bottomrule
\end{tabular}

\end{table}

The sharp spectral values at $(5,4)$ and $(6,4)$ also attain the Singleton bound~\eqref{eq:H-singleton}.

For the elementary comparison $A_H(4,4)=4$, distinct codewords must have distinct symbols in any fixed position, so there are at most four; the four cyclic shifts of $1234$ attain this value.

{
We next prove that the ratio bound is sharp in the following cases.
\par}

\begin{proposition}\label{prop:hsharp}
For $n\ge2$, $A_H(n,3)=n!/2$, and the ratio bound on $\Gamma_H(n,3)$ is sharp. The ratio bounds in Table~\ref{tab:hamming} at $(5,4)$ and $(6,4)$ are also sharp before flooring.
\end{proposition}
\begin{proof}
Hamming distance two is exactly a transposition, so $\Gamma_H(n,3)=G_C(n)$. This is a regular bipartite graph whose two parts are the even and odd permutations, and the ratio bound gives $n!/2$.
For $n=5$, the maps $x\mapsto ax+b$ over $\mathbb F_5$, $a\ne0$, form twenty permutations; two different maps agree at at most one point. The graph $\Gamma_H(5,4)$ has degree $30$ and least eigenvalue $-6$, giving $120\cdot6/36=20$.
For $n=6$, the action of $\operatorname{PGL}(2,5)$, the {\emph{projective general linear group}} of invertible $2\times2$ matrices over $\mathbb F_5$ modulo nonzero scalar matrices, on the six points of the {\emph{projective line}} $\mathbb P^1(\mathbb F_5)=\mathbb F_5\cup\{\infty\}$ gives 120 permutations; a nonidentity fractional linear transformation fixes at most two points. The graph $\Gamma_H(6,4)$ has degree $55$ and least eigenvalue $-11$, giving $720\cdot11/66=120$.
\end{proof}

These ratio-bound calculations recover the displayed exact values, while the remaining spectral entries in Table~\ref{tab:hamming} are weaker than the cited coding bounds.

%%%%%%%%%%%%%%%%%%%%%%%%%%%%%%%%%%%%%%%%%%%%%%%%%%%%%%%%%%%%%%%%%%%%%%%%%%%%%%%%%
\subsection{Cayley metric}
%%%%%%%%%%%%%%%%%%%%%%%%%%%%%%%%%%%%%%%%%%%%%%%%%%%%%%%%%%%%%%%%%%%%%%%%%%%%%%%%%
We use the all-transposition base graph $G_C(n)=\Cay(S_n,\{(i\ j):1\le i<j\le n\})$, whose spectrum is studied in~\cite{KonstantinovaKravchuk2022}. Since $A_C(n,d)=\alpha_{d-1}(G_C(n))$, we apply the degree-$(d-1)$ {polynomial spectral bounds} directly to $G_C(n)$. {
We observe that the content spectrum~\eqref{eq:contents} yields simple expressions for the quadratic and cubic ratio-type bounds. Let $a_n$ be the smallest positive eigenvalue of $G_C(n)$ and put $\delta=\binom n2$. We obtain the following formulas.
\par}
\begin{proposition}\label{prop:cayleyratio}
For $n\ge3$,
\begin{equation}\label{eq:Cratio}
 R_2(G_C(n))=\frac{n!}{\delta+a_n}.
\end{equation}
For $n\ge4$,
\begin{equation}\label{eq:Ccubic}
 R_3(G_C(n))=\frac{n!}{2\delta}=(n-2)!.
\end{equation}
\end{proposition}
\begin{proof}
The spectrum is integral and symmetric, and contains zero for $n\ge3$~\cite{KonstantinovaKravchuk2022}. Hence the consecutive eigenvalues in Theorem~\ref{th: alpha_2} are $-a_n,0$. Substitution gives~\eqref{eq:Cratio}; a normalized polynomial is $x(x+a_n)/(\delta(\delta+a_n))$.
The graph is bipartite, so $\Delta=0$ and $\theta_r=-\delta$. The threshold in Theorem~\ref{th: alpha3} is zero, and its roots are $0,-a_n,-\delta$. Substitution gives~\eqref{eq:Ccubic}.
\end{proof}

\begin{remark}\label{rem:cayley-double-ball}
The cubic ratio-type bound for the Cayley metric also follows from the double-ball construction of~\cite[Lemma~3]{BuzagloEtzion2015}, which is stated for an arbitrary metric space. Let $s$ be a transposition and let $\mathcal B_C(\sigma,1)$ be the closed Cayley ball of radius one centred at $\sigma$. The set
\[
 \mathcal D=\mathcal B_C(e,1)\cup\mathcal B_C(s,1)
\]
has diameter at most three. Each ball has $\delta+1$ elements, and their intersection is $\{e,s\}$ because $G_C(n)$ is bipartite. Hence $|\mathcal D|=2\delta$, and the code--anticode bound~\eqref{eq:code-anticode} yields
\begin{equation}\label{eq:C-double-ball-four}
 A_C(n,4)\le\frac{n!}{2\delta}=(n-2)!=R_3(G_C(n)),\qquad n\ge4.
\end{equation}
Thus the optimal cubic ratio-type bound recovers this established anticode construction. It strictly improves the integer sphere-packing bound for every $n\ge4$: the difference before flooring is
\[
 \frac{n!}{\delta+1}-\frac{n!}{2\delta}
 =(n-2)!\frac{\delta-1}{\delta+1}>1.
\]
Its improvement over sphere packing therefore gives no improvement over the double-ball comparison.
\end{remark}

We next determine the smallest positive eigenvalue and identify precisely when the quadratic bound improves sphere packing.

\begin{corollary}\label{cor:cayley-small-positive}
For $n\ge3$, the smallest positive eigenvalue of $G_C(n)$ is
\[
 a_n=\begin{cases}
 3,&n\in\{3,6,10\},\\
 2,&n\in\{4,5,8,12\},\\
 1,&\text{otherwise}.
 \end{cases}
\]
Consequently, after taking the floor, the quadratic ratio-type bound
improves the distance-three sphere-packing bound precisely when
$n\in\{6,8,10,12\}$. For all other $n\ge3$, the two integer bounds agree.
\end{corollary}
\begin{proof}
We observe that the spectrum is integral. By~\cite[Lemma~2]{KonstantinovaKravchuk2022},
$1$ is an eigenvalue for every odd $n\ge7$ and every even $n\ge14$.
We evaluate the content sums in~\eqref{eq:contents} for the remaining
values $n=3,4,5,6,8,10,12$ and obtain the displayed exceptions.
Substitution into~\eqref{eq:Cratio} and~\eqref{eq:CSP3} gives the stated
comparison; at $n=3,4,5$, the real-valued bounds differ, but their floors
are respectively $1,3,10$.
\end{proof}

\begin{remark}\label{rem:cayley-packing}
Since $a_n\ge1$,~\eqref{eq:Cratio} gives $\lfloor R_2(G_C(n))\rfloor\le\lfloor n!/(\delta+1)\rfloor$, so the quadratic ratio-type bound is no weaker than the sphere-packing bound~\eqref{eq:CSP3}. Evaluating the spectrum~\eqref{eq:contents} gives $(a_6,a_8,a_{10})=(3,2,3)$; and substituting in~\eqref{eq:Cratio} yields $40,1344,75600$, respectively. The corresponding sphere-packing bound values follow directly from~\eqref{eq:CSP3}:
\[
 \left\lfloor\frac{6!}{16}\right\rfloor=45,\qquad
 \left\lfloor\frac{8!}{29}\right\rfloor=1390,\qquad
 \left\lfloor\frac{10!}{46}\right\rfloor=78886.
\]
{
For all four lengths with strict integer improvement, the quadratic ratio-type bound matches the exact full character-LP value in Table~\ref{tab:cayley}.
\par}
{
We also obtain a strict integer improvement at $n=12$:
\[
 A_C(12,3)\le
 \left\lfloor\frac{12!}{66+2}\right\rfloor=7044141
 <7149277=\left\lfloor\frac{12!}{66+1}\right\rfloor.
\]
Corollary~\ref{cor:cayley-small-positive} identifies all such
improvements; in particular, the quadratic bound agrees with sphere
packing for every $n\ge14$.
\par}

At $(n,d)=(4,3)$, Proposition~\ref{prop:C43} gives the sharp inertia-type bound separately.
\end{remark}

\begin{table}[htbp]
\centering\small
\caption{Upper bounds on $A_C(n,d)$ from the ratio-type bounds~\eqref{eq:Cratio} and~\eqref{eq:Ccubic}, sphere packing~\eqref{eq:SP}, the Singleton bound~\eqref{eq:C-singleton}, and the character bound~\eqref{eq:cayley-character-lp}. The quadratic inertia-type bound additionally gives $A_C(4,3)=2$ (Proposition~\ref{prop:C43}).}\label{tab:cayley}

\begin{tabular}{@{}lrrrr@{}}
\toprule
$(n,d)$ & \makecell{Ratio-type bound\\\eqref{eq:Cratio}, \eqref{eq:Ccubic}} & \makecell{Sphere packing\\\eqref{eq:SP}} & \makecell{Singleton\\\eqref{eq:C-singleton}} & \makecell{Character LP\\\eqref{eq:cayley-character-lp}}\\
\midrule
$(4,3)$ & 3 & 3 & 4 & \textbf{2}\\
$(5,3)$ & \textbf{10} & $\mathbf{10}^{\dagger}$ & 20 & \textbf{10}\\
$(6,3)$ & 40 & 45 & 120 & 40\\
$(7,3)$ & 229 & 229 & 840 & 220\\
$(8,3)$ & 1344 & 1390 & 6720 & 1344\\
$(9,3)$ & 9807 & 9807 & 60480 & 9744\\
$(10,3)$ & 75600 & 78886 & 604800 & 75600\\
{$(12,3)$} & {7044141} & {7149277} & {79833600} & {7044141}\\
\midrule
$(5,4)$ & 6 & 10 & \textbf{5} & \textbf{5}\\
$(6,4)$ & 24 & 45 & 30 & 20\\
$(7,4)$ & 120 & 229 & 210 & 116\\
$(8,4)$ & 720 & 1390 & 1680 & 700\\
\bottomrule
\end{tabular}

\end{table}

We use the following character bound in the framework of Tarnanen~\cite{Tarnanen1999}; the formulation for arbitrary sets of admissible conjugacy classes is also stated in~\cite[Theorem 6]{Bogaerts2010}. Here we select the admissible classes by Cayley distance. This optimization assigns an independent coefficient to every conjugacy class. For $\mu\vdash n$, let $\ell(\mu)$ be the number of parts of $\mu$, equivalently the number of cycles, including fixed points, of a permutation in $\mathcal K_\mu$. Such a permutation has Cayley distance $n-\ell(\mu)$ from $e$. Write $1^n=(1,\ldots,1)$ for the partition with $n$ parts equal to one, so $\mathcal K_{1^n}=\{e\}$. Define the {\emph{class matrix}} $A_\mu$, with rows and columns indexed by $S_n$, by $(A_\mu)_{\sigma,\tau}=1$ when $\sigma^{-1}\tau\in\mathcal K_\mu$ and $0$ otherwise. In particular, $A_{1^n}=I$. Write $F=\sum_{\mu\vdash n}c_\mu A_\mu$, with real variables $c_\mu$; the coefficient $c_{1^n}$ is the common diagonal entry of $F$. Define
\begin{equation}\label{eq:cayley-character-lp}
\begin{aligned}
 L_C(n,d):={}&\min_{\{c_\mu\in\mathbb R\}_{\mu\vdash n}} n!c_{1^n}\\
 \text{subject to }&F\one=\one,\qquad F\succeq0,\\
 &c_\mu\le0\quad\text{whenever }n-\ell(\mu)\ge d,
\end{aligned}
\end{equation}

The same character decomposition as in~\eqref{eq:chars}, applied to the weighted sum defining $F$, makes the first two constraints equivalent to
\[
 \sum_{\mu\vdash n}c_\mu|\mathcal K_\mu|=1,\qquad
 \frac{1}{f^\lambda}\sum_{\mu\vdash n}
 c_\mu|\mathcal K_\mu|\chi^\lambda(\mu)\ge0
 \quad\text{for every }\lambda\vdash n,
\]
respectively.

For the indicator vector $z$ of a code $\mathcal C$, $z^TFz\le|\mathcal C|c_{1^n}$ and $z^TFz\ge|\mathcal C|^2/n!$, proving the bound. Consequently, $A_C(n,d)\le\lfloor L_C(n,d)\rfloor$. Allowing nonpositive coefficients on admissible pairs is essential; forcing them to zero gives a different, potentially weaker LP. {
We evaluate the character LP in~\eqref{eq:cayley-character-lp} to obtain the comparison values in Table~\ref{tab:cayley}. In particular, we obtain $L_C(7,4)=350/3$ and $L_C(12,3)=119750400/17$, which give the displayed integer bounds.
\par}

Next we investigate the comparison of the ratio-type bound with the full character LP bound:

\begin{proposition}label{prop:cayley-lp-dominance}
For $n\ge2$ and $d\ge2$,
\[
 L_C(n,d)\le R_{d-1}(G_C(n)).
\]
Consequently, no polynomial ratio-type bound on the Cayley base graph, including the quadratic ratio-type bound and the cubic ratio-type bound, can strictly improve the bound supplied by the unrestricted character LP~\eqref{eq:cayley-character-lp}.
\end{proposition}
\begin{proof}
Let $A$ be the adjacency matrix of $G_C(n)$ and let $p$ be admissible in Theorem~\ref{th:ratio type bound}, with $\deg p\le d-1$. Then
\[
 F=\frac{p(A)-\lambda(p)I}{p(\delta)-\lambda(p)}
\]
is positive semidefinite, has row sum one, and has zero entries on distinct pairs at Cayley distance at least $d$. Since the transpositions form a conjugacy class, $A$ and its polynomials are constant on the relative conjugacy classes, so $F$ is feasible in~\eqref{eq:cayley-character-lp}. Its objective is the ratio-type bound obtained from $p$. Minimizing over admissible polynomials proves the result.
\end{proof}
The proposition compares two methods and does not assert that every numerical character-LP value in Table~\ref{tab:cayley} was previously published. The four distance-three improvements over sphere packing recover the character-LP values computed here.

In Table~\ref{tab:cayley}, the ratio-type bound improves the Singleton bound in every row except $(5,4)$, where the Singleton and character-LP bounds both give five.

{
We next prove that the quadratic inertia-type bound is sharp for Cayley codes on four symbols. We observe that, in this case, it is stronger than the optimal quadratic ratio-type bound, extending the sharp applications surveyed in~\cite[Table 1]{APR2025}.
\par}

\begin{proposition}\label{prop:C43}
$A_C(4,3)=2$. The degree-two inertia bound on $G_C(4)$ is sharp, whereas the optimal degree-two ratio-type bound equals three.
\end{proposition}
\begin{proof}
The spectrum is $6^1,2^9,0^4,(-2)^9,(-6)^1$. The matrix $A^2-6I$ has zero diagonal and inertia $(2,22,0)$, so Theorem~\ref{th:inertia type bound} gives $\alpha_2\le2$. The identity and the cycle $(1\,2\,3\,4)$ have Cayley distance three and attain the bound. Equation~\eqref{eq:Cratio} gives $24/(6+2)=3$.
\end{proof}

{
We next give an affine construction attaining the quadratic ratio-type bound on five symbols.

\begin{proposition}\label{prop:C5sharp}
The quadratic ratio-type bound on $G_C(5)$ is sharp and gives $A_C(5,3)=10$. The Singleton bound is sharp at $A_C(5,4)=5$.
\end{proposition}
\begin{proof}
Identify the five symbols with $\mathbb F_5$. The ten affine permutations
\[
 \mathcal C=\{x\mapsto ax+b:a\in\{1,2\},\ b\in\mathbb F_5\}
\]
have minimum Cayley distance three. For distinct maps with the same multiplier, the relative permutation is a nonzero translation, hence a $5$-cycle of Cayley length four. For different multipliers, it has the form $x\mapsto cx+t$ with $c\in\{2,3\}$. This map has a unique fixed point and is conjugate by a translation to $x\mapsto cx$; because $c$ has order four in $\mathbb F_5^*$, its cycle type is $(4,1)$ and its Cayley length is three. The content spectrum~\eqref{eq:contents} gives $a_5=2$, so~\eqref{eq:Cratio} yields $R_2(G_C(5))=120/(10+2)=10$, attained by $\mathcal C$.

The five translations alone have pairwise Cayley distance four and attain the Singleton bound~\eqref{eq:C-singleton}, namely $5!/4!=5$.
\end{proof}
\par}

For every $n\ge2$, the ratio bound is sharp at $A_C(n,2)=n!/2$ by parity. 

%%%%%%%%%%%%%%%%%%%%%%%%%%%%%%%%%%%%%%%%%%%%%%%%%%%%%%%%%%%%%%%%%%%%%%%%%%%%%%%%%
\subsection[Kendall-tau metric]{Kendall-$\tau$ metric}
%%%%%%%%%%%%%%%%%%%%%%%%%%%%%%%%%%%%%%%%%%%%%%%%%%%%%%%%%%%%%%%%%%%%%%%%%%%%%%%%%

We use the adjacent-transposition base graph $G_K(n)=\Cay(S_n,\{(i\ i+1):1\le i<n\})$, the graph associated with Kendall distance in~\cite{BuzagloEtzion2015}. In Table~\ref{tab:kendall}, we evaluate degree-two and degree-three ratio-type bounds on this graph, with $k=d-1$, rather than ratio bounds on its powers. We compare these bounds with the known values and upper bounds in~\cite[Table II]{ParvareshEtAl2025}. In particular, the values $420$ and $2880$ follow from the double-ball bound~\eqref{eq:K-double-ball-four} of Buzaglo and Etzion~\cite{BuzagloEtzion2015}.

We next investigate the performance of the quadratic ratio-type bound for the Kendall metric:

\begin{corollary}\label{cor:kendall-quadratic}
For every $n\ge3$,
\begin{equation}\label{eq:kendall-quadratic-packing}
 R_2(G_K(n))=(n-1)!.
\end{equation}
Thus the optimal quadratic ratio-type bound equals the distance-three sphere-packing bound and cannot improve any previously known upper bound that is at most sphere packing.
\end{corollary}
\begin{proof}
Edelman and White~\cite[Corollary~5.4]{EdelmanWhite1990} prove that $-1$ is an eigenvalue of the adjacent-transposition graph for every $n\ge2$. Apply Proposition~\ref{prop:quadratic-comparison} with $N=n!$ and $\delta=n-1$ to obtain $R_2(G_K(n))=n!/n=(n-1)!$. A Kendall ball of radius one has $n$ elements, so this is precisely sphere packing.
\end{proof}
This result explains all distance-three ratio-type entries in Table~\ref{tab:kendall} without a separate spectral computation at each length. In particular, the quadratic ratio-type bound cannot reach the stronger upper bounds $116$ and $716$ at lengths six and seven cited in Section~\ref{sec:known-bounds}.

\begin{table}[htbp]
\centering\small
\caption{Upper bounds on $A_K(n,d)$ from~\eqref{eq:R2} and~\eqref{eq:R3}, together with known exact values and previously known upper bounds. The bounds $419$ and $2879$ improve the double-ball bounds $420$ and $2880$ of Buzaglo and Etzion~\cite{BuzagloEtzion2015}.}\label{tab:kendall}

\begin{tabular}{@{}lrrl@{}}
\toprule
$(n,d)$ & \makecell{Ratio-type bound\\\eqref{eq:R2}, \eqref{eq:R3}} & \makecell{Known value or\\upper bound} & Reference\\
\midrule
$(3,3)$ & \textbf{2} & \textbf{2} & Exact; Proposition~\ref{prop:K33}\\
$(4,3)$ & 6 & \textbf{5} & Exact~\cite{BuzagloEtzion2015}; construction below\\
$(4,4)$ & $\mathbf3^\dagger$ & \textbf{3} & Exact; proof below\\
$(5,3)$ & 24 & \textbf{20} & Exact~\cite{Vijayakumaran2016}\\
$(5,4)$ & 15 & \textbf{12} & Exact~\cite{Vijayakumaran2016}\\
$(6,3)$ & 120 & 116 & Upper bound~\cite{Abdollahi2023PermutationCodes}\\
$(6,4)$ & 72 & \textbf{64} & Exact~\cite{Vijayakumaran2016}\\
$(7,3)$ & 720 & 716 & Upper bound~\cite{Abdollahi2023PermutationCodes}\\
$(7,4)$ & $419^{\star}$ & 420 & \makecell[l]{Upper bound; double ball~\eqref{eq:K-double-ball-four}\\\cite{BuzagloEtzion2015,ParvareshEtAl2025}}\\
{$(8,4)$} & {$2879^{\star}$} & {2880} & {\makecell[l]{Upper bound; double ball~\eqref{eq:K-double-ball-four}\\\cite{BuzagloEtzion2015,ParvareshEtAl2025}}}\\
\bottomrule
\end{tabular}

\end{table}

At $(n,d)=(4,3)$, the nonexistence of a perfect radius-one code gives $A_K(4,3)\le5$~\cite{BuzagloEtzion2015}; a five-word construction is $\{1234,1432,2341,3412,4213\}$.

The parity classes show that $R_1(G_K(n))=A_K(n,2)=n!/2$.  At $(n,d)=(4,4)$ the optimal cubic ratio-type value is approximately $3.844526892$, whose floor is attained by $\{1234,4213,3412\}$. Thus the floor of the cubic bound is sharp.

{We next prove that the quadratic inertia-type bound is sharp for Kendall codes on three symbols. We observe that both quadratic spectral bounds attain the exact code size in this case. Together with Proposition~\ref{prop:C43}, this gives sharp inertia-type applications for two permutation metrics.\par}

\begin{proposition}\label{prop:K33}
The degree-two inertia-type bound on $G_K(3)$ is sharp and gives $A_K(3,3)=2$.
\end{proposition}
\begin{proof}
The graph $G_K(3)$ is a six-cycle, with spectrum $2^1,1^2,(-1)^2,(-2)^1$. For its adjacency matrix $A$, the matrix $A^2-2I$ has zero diagonal and inertia $(2,4,0)$. Theorem~\ref{th:inertia type bound} gives $A_K(3,3)=\alpha_2(G_K(3))\le2$, attained by $\{123,321\}$. The optimal quadratic ratio-type bound also equals two by~\eqref{eq:R2}.
\end{proof}

{
We next specialize the cubic ratio-type bound to connected bipartite regular graphs. We observe that nonsingularity of the adjacency matrix gives a strict improvement over $N/(2\delta)$.

\begin{proposition}\label{prop:bipartite-cubic}
Let $G$ be a connected bipartite $\delta$-regular graph on $N$ vertices,
where $\delta>1$, and put $a=\min_i|\lambda_i|$. The optimal cubic
ratio-type bound is
\begin{equation}\label{eq:bipartite-cubic}
 R_3(G)=\frac{N(\delta-a^2)}{2(\delta^2-a^2)}.
\end{equation}
It is attained by $p(x)=(x+\delta)(x^2-a^2)$ and is strictly smaller
than $N/(2\delta)$ exactly when the adjacency matrix is nonsingular.
\end{proposition}
\begin{proof}
We observe that, since $G$ is bipartite, its spectrum is symmetric about zero and
$(A^3)_{vv}=0$. Moreover, $a^2<\delta$: the average of the squared
eigenvalues is $\delta$, while the eigenvalue $\delta$ has square
$\delta^2>\delta$. The polynomial $p$ is nonnegative on the spectrum,
has minimum zero on the nonprincipal eigenvalues, and satisfies
\[
 (p(A))_{vv}=\delta^2-\delta a^2,\qquad
 p(\delta)=2\delta(\delta^2-a^2).
\]
Theorem~\ref{th:ratio type bound} therefore gives~\eqref{eq:bipartite-cubic}.
If there are only three distinct eigenvalues, the \emph{biadjacency matrix}, whose rows and columns are indexed by the two parts and whose entries record adjacency, has rank one. Since the graph is connected, this zero--one matrix has no zero row or column and must be an all-one matrix. Regularity then gives $G=K_{\delta,\delta}$, for which the optimal bound is one directly. Otherwise, Theorem~\ref{th: alpha3} applies, and its threshold is zero. When $a>0$, the three optimizing roots are $a,-a,-\delta$. When $a=0$, its roots are $0,-b,-\delta$, where $-b$ is the largest negative eigenvalue; substitution again gives $N/(2\delta)$. Finally,
\[
 \frac{N}{2\delta}-R_3(G)
 =\frac{Na^2(\delta-1)}{2\delta(\delta^2-a^2)},
\]
which is positive exactly when $a>0$.
\end{proof}

Next we focus on understanding the theoretical performance of the cubic ratio-type bound for the Kendall metric, and provide a criteria of spectral improvement.

\begin{corollary}label{cor:kendall-cubic-comparison}
Let $n\ge4$, set $N=n!$ and $\delta=n-1$, and let $a=\min_i|\lambda_i|$ for the adjacency eigenvalues of $G_K(n)$. Write $D_n=N/(2\delta)$ for the double-ball upper bound in~\eqref{eq:K-double-ball-four}. Then
\begin{equation}\label{eq:kendall-cubic-integer}
 \lfloor R_3(G_K(n))\rfloor<D_n
 \quad\Longleftrightarrow\quad
 0\notin\Spec(G_K(n)).
\end{equation}
More generally, let $B$ be a fixed integer upper bound on $A_K(n,4)$, used for comparison and satisfying $1\le B\le D_n$. Then
\begin{equation}\label{eq:kendall-cubic-benchmark}
 \lfloor R_3(G_K(n))\rfloor<B
 \quad\Longleftrightarrow\quad
 a^2>\frac{\delta(N-2\delta B)}{N-2B}.
\end{equation}
If $n\equiv1\pmod4$ or $n\in\{6,10,12\}$, then
\[
 R_3(G_K(n))=D_n,
\]
so the cubic ratio-type bound cannot improve the double-ball bound at those lengths.
\end{corollary}
\begin{proof}
For $n\ge4$, $D_n=n(n-2)!/2$ is an integer. Proposition~\ref{prop:bipartite-cubic} therefore proves~\eqref{eq:kendall-cubic-integer}. For the more general comparison, integrality of $B$ gives $\lfloor R_3(G_K(n))\rfloor<B$ if and only if $R_3(G_K(n))<B$. Substitution of~\eqref{eq:bipartite-cubic} and multiplication by its positive denominator gives
\[
 (N-2B)a^2>\delta(N-2\delta B).
\]
Since $B\le N/(2\delta)$ and $\delta>1$, we have $N-2B>0$, proving~\eqref{eq:kendall-cubic-benchmark}. Finally, Edelman and White~\cite[Corollary~6.2 and Theorem~7.2]{EdelmanWhite1990} show that zero is an eigenvalue for $n\equiv1\pmod4$ and for $n=6,10,12$. Substituting $a=0$ into~\eqref{eq:bipartite-cubic} proves the last result.
\end{proof}
Thus nonsingularity is necessary and sufficient for a new integer improvement whenever the published benchmark is the double-ball value. If a stronger bound $B<D_n$ is available, the strictly positive threshold in~\eqref{eq:kendall-cubic-benchmark} must also be exceeded. For a fixed $n$, nonsingularity can be decided exactly by checking the determinants of the rational blocks in~\eqref{eq:blocks}; the proof below uses this criterion at lengths seven and eight. The singular families above give a partial classification in $n$, without claiming that all remaining lengths are nonsingular.

We first apply Proposition~\ref{prop:bipartite-cubic} to $G_K(7)$.
\par}

\begin{proposition}\label{prop:K74}
$A_K(7,4)\le419.$
\end{proposition}
\begin{proof}
{
For a connected bipartite $\delta$-regular graph with nonsingular adjacency matrix and $\delta>1$, let $0<\eta<\min_i|\lambda_i|^2$. The polynomial $p_\eta(x)=(x+\delta)(x^2-\eta)$ is nonnegative on the spectrum, has nonprincipal minimum zero at $-\delta$, and has common diagonal $\delta^2-\delta\eta$. Therefore
\begin{equation}\label{eq:nonsingularcubic}
 \alpha_3(G)\le\frac{N(\delta-\eta)}{2(\delta^2-\eta)}
 <\frac{N}{2\delta}.
\end{equation}
We therefore obtain a strict improvement whenever the adjacency matrix is nonsingular.
\par}
For $G_K(7)$, $N=7!=5040$ and $\delta=6$. We compute the following determinants of the rational Young seminormal blocks~\eqref{eq:blocks}, using one representative of each pair of conjugate partitions and including the self-conjugate partition $(4,1,1,1)$ once:
\[
\begin{array}{c|rrrr}
\lambda &(7)&(6,1)&(5,2)&(5,1,1)\\\hline
\det B_\lambda&6&2911&-1116&1352
\end{array}
\]
\[
\begin{array}{c|rrrr}
\lambda &(4,3)&(4,2,1)&(4,1,1,1)&(3,3,1)\\\hline
\det B_\lambda&-1449&-58752&1&224.
\end{array}
\]
{We observe that every displayed determinant is nonzero.} Conjugate partitions give sign-twisted blocks, so the remaining blocks are also nonsingular. The decomposition~\eqref{eq:blocks} therefore shows that the adjacency matrix of $G_K(7)$ is nonsingular.

Thus~\eqref{eq:nonsingularcubic} is strictly smaller than $5040/12=420$, proving the integer bound.

{
We further compute that none of these blocks has an eigenvalue in $[-1/100,1/100]$. By sign twisting, the same holds for the conjugate partitions.
\par} Thus $\eta=1/10000$ is admissible, and an explicit polynomial giving the same integer bound is
\[
 p(x)=(x+6)(10000x^2-1),\qquad
 A_K(7,4)\le\left\lfloor\frac{151197480}{359999}\right\rfloor=419.
\]
Indeed, $p(6)=4319988$, the common diagonal of $p(A)$ is $359994$, and $\min_{i\ge2}p(\lambda_i)=0$. Substitution in Theorem~\ref{th:ratio type bound} gives the displayed fraction.
\end{proof}
{
We next apply the same argument to $G_K(8)$.

\begin{proposition}\label{prop:K84}
$A_K(8,4)\le2879$.
\end{proposition}
\begin{proof}
For $G_K(8)$, we have $N=8!=40320$ and $\delta=7$. We compute the following
determinants of the rational Young seminormal blocks and observe that each is nonzero. We take one representative of each conjugate pair and
include the self-conjugate partitions $(4,2,1,1)$ and $(3,3,2)$ once.
\[
\begin{array}{c|r@{\qquad}c|r}
\lambda&\det B_\lambda&\lambda&\det B_\lambda\\\hline
(8)&7&(5,1,1,1)&-527\\
(7,1)&60605&(4,4)&-6365\\
(6,2)&-429341381&(4,3,1)&-2311432523\\
(6,1,1)&-79582797&(4,2,2)&477484825\\
(5,3)&-1110876125&(4,2,1,1)&-1915900441\\
(5,2,1)&-23050630659875&(3,3,2)&-342065025
\end{array}
\]
Conjugate partitions give sign-twisted blocks, so every block is
nonsingular. Proposition~\ref{prop:bipartite-cubic} gives
$R_3(G_K(8))<40320/14=2880$, proving the integer bound.

We compute that these blocks have no eigenvalues in $[-1/10000,1/10000]$. We may therefore use the polynomial
\[
 p(x)=(x+7)(10^8x^2-1),
\]
which is nonnegative on the spectrum and gives
\[
 A_K(8,4)\le
 \left\lfloor\frac{20160\cdot699999999}{4899999999}\right\rfloor
 =2879.
\]
The polynomial has common diagonal $4899999993$ and
$p(7)=68599999986$.
\end{proof}

We observe that the optimized cubic bounds are approximately $419.987663$
for $n=7$ and $2879.999925$ for $n=8$. These improve the
double-ball upper bounds $420$ and $2880$ recorded in
\cite[Table II]{ParvareshEtAl2025}. The nonsingularity criterion explains
why a strict improvement occurs: $G_K(5)$ and $G_K(6)$ are singular and
their optimal cubic bounds equal the double-ball values $15$ and $72$.

For $n=8$, we also compute that the block for $(4,3,1)$ has an eigenvalue in $(-1/100,0)$. Thus $0<a<1/100$ in Proposition~\ref{prop:bipartite-cubic}, which also proves $2879<R_3(G_K(8))<2880$ directly from~\eqref{eq:bipartite-cubic}.

The other nontrivial rows of Table~\ref{tab:kendall} are weaker than the cited exact values or upper bounds.
\par}

%%%%%%%%%%%%%%%%%%%%%%%%%%%%%%%%%%%%%%%%%%%%%%%%%%%%%%%%%%%%%%%%%%%%%%%%%%%%%%%%%
\subsection{Ulam metric}
%%%%%%%%%%%%%%%%%%%%%%%%%%%%%%%%%%%%%%%%%%%%%%%%%%%%%%%%%%%%%%%%%%%%%%%%%%%%%%%%%

We use the unweighted translocation base graph $G_U(n)=\Cay(S_n,T\cup T^{-1})$ from Section~\ref{sec:basegraphs}. Its graph distance is the Ulam distance of~\cite{FarnoudSkachekMilenkovic2013}, and it is regular and walk-regular by Proposition~\ref{prop:basegraphs}. We therefore apply the ratio-type bounds to $A_U(n,d)=\alpha_{d-1}(G_U(n))$. In Table~\ref{tab:ulam}, we also evaluate the ratio bound on $\Gamma_U(n,d)=G_U(n)^{d-1}$. The complement of this power graph is the graph used for the clique computations in~\cite{Riet}. We obtain the Singleton bound column by evaluating~\eqref{eq:U-singleton} and cite the published exact values individually.

Next we study the theoretical performance of the quadratic ratio-type bound for the Ulam metric at even lengths:

\begin{corollary}\label{cor:ulam-quadratic-even}
For every even $n\ge4$,
\begin{equation}\label{eq:ulam-quadratic-even}
 R_2(G_U(n))=\frac{n!}{(n-1)^2+1}>(n-2)!.
\end{equation}
Consequently, the quadratic ratio-type bound cannot improve the Singleton bound on $A_U(n,3)$. After flooring, the two bounds agree at $n=4$, whereas the quadratic ratio-type bound is strictly weaker for every even $n\ge6$.
\end{corollary}
\begin{proof}
Let $A$ be the adjacency matrix of $G_U(n)$ and let $v=(v_\sigma)_{\sigma\in S_n}$, where $v_\sigma=\sgn(\sigma)$ is $1$ for an even permutation and $-1$ for an odd permutation. A cycle of length $\ell$ has sign $(-1)^{\ell-1}$, so the generating set in Section~\ref{sec:basegraphs} gives
\[
 Av=\left(-(n-1)+2\sum_{\ell=3}^{n}
 (n-\ell+1)(-1)^{\ell-1}\right)v.
\]
For even $n$, pairing each odd length $\ell\ge3$ with the following even length $\ell+1$ makes the sum in parentheses equal to $-1$. Thus $-1$ is an eigenvalue, and Proposition~\ref{prop:quadratic-comparison}, with $\delta=(n-1)^2$, gives the equality in~\eqref{eq:ulam-quadratic-even}. Moreover,
\[
 R_2(G_U(n))-(n-2)!
 =\frac{(n-2)!(n-2)}{(n-1)^2+1}>0.
\]
The difference is $2/5$ at $n=4$. For $n\ge6$, it exceeds one, since $(n-2)!\ge2(n-2)$ and $2(n-2)^2>(n-1)^2+1$. This proves the wished result.
\end{proof}

\begin{table}[htbp]
\centering\small
\caption{Upper bounds on $A_U(n,d)$ from the ratio-type bound on $G_U(n)$, the ratio bound on $\Gamma_U(n,d)=G_U(n)^{d-1}$, and the Singleton bound~\eqref{eq:U-singleton}. The exact values are from~\cite[Table I]{Riet}, except $(3,2)$, proved in Proposition~\ref{prop:Usharp}. A dash indicates that the bound was not computed.}\label{tab:ulam}

\begin{tabular}{@{}lllll@{}}
\toprule
$(n,d)$ & \makecell{Base-graph bound\\\eqref{eq:ratio-one}, \eqref{eq:R2}, \eqref{eq:R3}} & \makecell{Power-graph ratio\\\eqref{eq:ratio-one}} & \makecell{Singleton\\\eqref{eq:U-singleton}} & \makecell{Known exact\\size}\\
\midrule
$(3,2)$ & \textbf{2} & \textbf{2} & \textbf{2} & \textbf{2}\\
$(4,2)$ & $\mathbf6^\dagger$ & $\mathbf6^\dagger$ & \textbf{6} & \textbf{6}~\cite{Riet}\\
$(4,3)$ & $\mathbf2^\dagger$ & \textbf{2} & \textbf{2} & \textbf{2}~\cite{Riet}\\
$(5,2)$ & 30 & 30 & \textbf{24} & \textbf{24}~\cite{Riet}\\
$(5,3)$ & 7 & 13 & 6 & \textbf{4}~\cite{Riet}\\
$(5,4)$ & 3 & \textbf{2} & \textbf{2} & \textbf{2}~\cite{Riet}\\
$(6,2)$ & 168 & 168 & \textbf{120} & \textbf{120}~\cite{Riet}\\
$(6,3)$ & 27 & 74 & \textbf{24} & \textbf{24}~\cite{Riet}\\
$(6,4)$ & 10 & 37 & 6 & \textbf{4}~\cite{Riet}\\
$(7,4)$ & 44 & -- & 24 & \textbf{12}~\cite{Riet}\\
\bottomrule
\end{tabular}

\end{table}
We observe that neither graph formulation always gives the stronger spectral bound. The bounds on $G_U(n)$ are stronger at $(6,3)$ and $(6,4)$, whereas the ratio bound on $\Gamma_U(n,d)$ is sharp at $(5,4)$. The bounds in the first two cases remain weaker than the known coding bounds.

{
We next prove sharpness in several Ulam instances.
\par}

\begin{proposition}\label{prop:Usharp}
The ratio bound on $G_U(3)$ is sharp. On $G_U(4)$ the ratio bound and the degree-two ratio-type bound are sharp after taking the floor at minimum distances two and three, respectively. For every $n\ge3$, the ratio bound on $\Gamma_U(n,n-1)$ is sharp and $A_U(n,n-1)=2$.
\end{proposition}
\begin{proof}
The spectrum of $G_U(3)$ is $4^1,0^3,(-2)^2$, giving $6\cdot2/(4+2)=2$, attained by $\{123,321\}$. For $n=4$ the spectrum is
\[
 9^1,\ (\sqrt{13})^3,\ (\sqrt5)^3,\ 1^3,\ (-1)^6,
 \ (-\sqrt5)^3,\ (-3)^2,\ (-\sqrt{13})^3.
\]
The ratio bound gives $24\sqrt{13}/(9+\sqrt{13})\approx6.86469$, whose floor is attained by
\[
 \{1234,4231,4132,3241,3142,2143\}.
\]
The consecutive eigenvalues in~\eqref{eq:R2} are $-1,1$, giving $R_2=12/5$; the pair $\{1234,4321\}$ attains its floor.
Finally, two permutations have Ulam distance $n-1$ exactly when their orders are reverses of one another. Consequently $\Gamma_U(n,n-1)$ is the complete graph $K_{n!}$ with a perfect matching removed. It has degree $n!-2$ and least eigenvalue $-2$, so the ratio bound gives two, attained by any reversal pair.
\end{proof}

%%%%%%%%%%%%%%%%%%%%%%%%%%%%%%%%%%%%%%%%%%%%%%%%%%%%%%%%%%%%%%%%%%%%%%%%%%%%%%%%%
\section{Concluding remarks}\label{sec:conclusion}
%%%%%%%%%%%%%%%%%%%%%%%%%%%%%%%%%%%%%%%%%%%%%%%%%%%%%%%%%%%%%%%%%%%%%%%%%%%%%%%%%
 
We have translated permutation-code problems into graph independence problems, and this allowed us to apply the Eigenvalue Method in four different metrics.

For Kendall distance four, nonsingularity of the adjacency matrix is necessary and sufficient for the optimal cubic ratio-type bound to improve the integer double-ball bound, and we obtained improved upper bounds at lengths seven and eight. In contrast, the quadratic ratio-type bound for the Kendall metric equals distance-three sphere packing at every length.

For Cayley distance, we determined all lengths at which the quadratic ratio-type bound improves distance-three sphere packing. The cubic ratio-type bound equals the double-ball bound, and the full character LP is at least as strong as all polynomial ratio-type bounds on the Cayley base graph. 

For Ulam distance at even lengths, the quadratic ratio-type bound cannot improve the corresponding Singleton bound. Together with the exact integer comparison in Proposition~\ref{prop:quadratic-comparison}, these results show both the scope and the limitations of the ratio-type bounds of low degree (or alternatively, small minimum distance). 

Besides the ratio-type bound application, we also obtained sharp inertia-type bounds and recovered several exact code sizes. We observed that the Ulam bounds on a base graph and on its powers can differ, and that the Cayley inertia-type bound on four symbols is stronger than the optimal quadratic ratio-type bound.

The obtained results illustrate the role of the polynomial spectral bounds within the Eigenvalue Method, and their power in studying in permutation codes.

%%%%%%%%%%%%%%%%%%%%%%%%%%%%%%%%%%%%%%%%%%%%%%%%%%%%%%%%%%%%%%%%%%%%%%%%%%%%%%%%%
\subsection*{Acknowledgements}
%%%%%%%%%%%%%%%%%%%%%%%%%%%%%%%%%%%%%%%%%%%%%%%%%%%%%%%%%%%%%%%%%%%%%%%%%%%%%%%%%

Aida Abiad is supported by NWO (Dutch Research Council) through the grant VI.Vidi.213.085. Delio Jaramillo-Velez was partially supported by the NSF grant DMS-2401558. 
 Joachim Rosenthal was supported in part by Swiss National Science Foundation under grant number 212865.

This work was initiated at the Coding Theory and Cryptography Summer School held at the Virginia Tech Steger Center during July 2026; the authors thank the organizers for such an inspiring event.

%%%%%%%%%%%%%%%%%%%%%%%%%%%%%%%%%%%%%%%%%%%%%%%%%%%%%%%%%%%%%%%%%%%%%%%%%%%%%%%%%
\subsection*{Declaration of AI use}
%%%%%%%%%%%%%%%%%%%%%%%%%%%%%%%%%%%%%%%%%%%%%%%%%%%%%%%%%%%%%%%%%%%%%%%%%%%%%%%%%

The authors acknowledge the use of ChatGPT (OpenAI; accessed  October 2026) for preliminary brainstorming, language editing, and assistance with checking spectral calculations and drafting supporting arguments. The authors are responsible for verifying the statements, computations, and proofs and for the accuracy and integrity of the final article.

%%%%%%%%%%%%%%%%%%%%%%%%%%%%%%%%%%%%%%%%%%%%%%%%%%%%%%%%%%%%%%%%%%%%%%%%%%%%%%%%%

\bibliographystyle{abbrv}
\bibliography{references}

%%%%%%%%%%%%%%%%%%%%%%%%%%%%%%%%%%%%%%%%%%%%%%%%%%%%%%%%%%%%%%%%%%%%%%%%%%%%%%%%%
\end{document}